\documentclass{amsart}
\usepackage{amsmath}
\usepackage{amsfonts}
\usepackage{amssymb}
\usepackage{graphicx,xcolor}
\newtheorem{Theorem}{Theorem}[section]

\newtheorem{Lemma}[Theorem]{Lemma}
\newtheorem{Proposition}[Theorem]{Proposition}

\newtheorem{Definition}[Theorem]{Definition}
\newtheorem{rem}[Theorem]{Remark}
\newtheorem{example}[Theorem]{Example}

\newcommand{\R}{\mathbb R}

\newcommand{\N}{\mathbb N}
\newcommand{\Z}{\mathbb Z}

\newcommand{\C}{\mathcal{C}}
\newcommand{\F}{\mathcal{F}}
\newcommand{\p}{\mathcal{P}}

\newcommand{\D}{\mathcal{D}}

\newcommand{\rel}{\mathcal{R}}

\newcommand{\mcc}{\mathcal{C}}

\title[differential geometry of numerical schemes and weak solutions]{On the differential geometry of numerical schemes and weak solutions of functional equations.}
\author{Jean-Pierre Magnot}
\address{LAREMA, Universit\'e d'Angers, 2 Bd Lavoisier, 49045 Angers cedex 1, ´
France and Lyc\'ee Jeanne d'Arc \\ Avenue de Grande Bretagne \\ F-63000 Clermont-Ferrand}
\email{jean-pierr.magnot@ac-clermont.fr}
\begin{document}

\maketitle

\begin{abstract}
We exhibit  differential geometric structures that arise in numerical methods, based on the construction of Cauchy sequences, that are currently used to prove explicitly the existence of weak solutions to functional equations. We describe the geometric framework, highlight several examples and describe how two well-known proofs  fit with our setting. The first one is a re-interpretation of the classical proof of an implicit functions theorem in an ILB setting, for which our setting enables us to state an implicit functions theorem without additional norm estimates, and the second one is the finite element method of {{}a} Dirichlet problem where the set of triangulations appear as a smooth set of parameters. In both case, smooth dependence on the set of parameters is established. Before that, we develop the necessary theoretical tools, namely the notion of Cauchy diffeology on spaces of Cauchy sequences and a new generalization of the notion of tangent space to a diffeological space. 
\end{abstract}

\vskip 12pt
\textit{Keywords:} Diffeology, Cauchy sequences, functional equations, weak solutions, triangulations.

\textit{MSC (2010):} 46T20; 58C15; 58B10; 47J25
\section*{Introduction}
The aim of this paper is to describe some differential geometric properties of the analysis of weak solutions of functional equations, especially partial differential equations. We highlight a generalized differential geometric structure, called diffeology (we call it generalized differential geometry because this setting does not {{}involve} atlases), as well as a definition of (abstract) symmetries, based on the notion of numerical schemes that are commonly used in contructing explicit solutions to functional equations of the type $$F(u,q)=0$$ where $u$ is a function and $q$ is a parameter. We refine for this goal two theoretical aspects of diffeologies, that carry a language which is quite user-friendly. First we define a new tangent space on a diffeological space, and secondly we define a diffeology on Cauchy sequences, that we call Cauchy diffeology, that seems adapted to this setting. We illustrate the results and the settings of this paper by several examples. We finish with two worked-out examples adapted from well-known frameworks. First we analyze the problem of implicit functions under the light of weak solutions. This enables us to state it in the ILB setting \cite{Om}, for degree 0 map which do not carry additional norm estimates as in classical statements \cite{Ham,HN1971,Om}, {{}finishing the study initiated in \cite{Ma2020-1}.} The difference {{}between} our result with these classical appraoches is discussed in details, as well as the correspondence with the geometry of weak solutions. Secondly we complete the study of {{}the} degree 1 finite elements method for the Dirichelet problem. Here existence and uniqueness of solutions is well-known but we concentrate our efforts on the forgotten aspect of the diffeology of the space of triangulations of the domain, and we show how triangulations can be considered as a smooth space of parameters for numerical schemes built up from the finite elements method. 

Let us describe with more details the contents of this paper. We recall the necessary material on diffeologies in section 1.

In section \ref{T2}, we describe a diffeology, new to our knowledge, which appears as a refined diffeology of a diffeological space. This refined diffeology is inherited from the group of diffeomorphisms, and gives rise to a new (fourth) definition of tangent space of a diffeological space. These definitions are useful for the generalization of the notion of symmetries proposed in section \ref{3}. 

In section \ref{2}, we describe the so-called Cauchy diffeology on Cauchy sequences, for which 
\begin{enumerate}
	\item the limit map is smooth 
	\item the index maps $ev_k:(u_n)_{n \in \N} \mapsto u_k$ are smooth $(k \in \N).$
\end{enumerate}
We show that the two conditions are necessary, by well-chosen examples. {{} The good conditions for fitting with the technical requirements in the considered applications are found when the two conditions are gathered simultatenouly.} 

In section \ref{3}, we find another motivation for the diffeology: when solving numerically a PDE, we build a sequence $(u_n)$ which converges to the solution $u,$ but when $u_n$ is "close enough" to $u,$ computer representations of $u$ uses the approximate solution $u_n.$ Thus, {{} both} $u$ and $u_n$ need to be smooth {{} with respect to} the parameters and the initial conditions. After describing what can be a general setting for numerical methods {{} and for their symmetry groups}, we  address an open question on {{}paradoxical} solutions of Euler equations for perfect fluids {{}\cite{dLS1,dLS2,dLS3,Sch,Sh1,Sh2,Vil}}. Cauchy sequences appear mostly where convergence of sequences is needed. This tool is basically topological, as well as the notion of convergence. But especially in analysis of ordinary or partial differential equations, once Cauchy sequences {{}are built up to approximate} solutions, a classical question is the smooth dependence on {{} initial conditions and/or  parameters}. Out of a well-established manifold structure on the parameters or  initial conditions, and even if these structures are given, diffeologies appear as an easy way to formalize smoothness {{} with less technical constraints than in} the more rigid framework of manifolds. 

In section \ref{4}, we show that part of the classical hypothesis of classical (smooth) implicit functions theorems can be relaxed {{} by introducing} the diffeology of Cauchy sequences {{} in the topological approach initiated in \cite{Ma2020-1} where the choice to reduce the domain $D$ of the implicit function is made} instead of starting with strong estimates on the considered functions {{} in order to control better the nature of $D.$} We concentrate on an implicit function theorem on ILB sequences of Banach spaces $(E_i)$ and $(F_i)$ and smooth maps $f_i: O_i \subset E_i \times F_i \rightarrow F_i.$ {{} As in \cite{Ma2020-1},} uniform estimates on the family $f_i$ are not necessary to define a diffeological domain $D\subset \bigcap_{i \in \N } O_i$ and a smooth map $u : D \rightarrow \bigcap_{i \in \N}F_i$ such that $\forall i \in \N, f_i(x,u(x))=0.$
  We finish with a corresponding ``free of estimates'' Fr\"obenius theorem, again with the only help of rewritten classical proofs.  

In section \ref{5}, we show how one of the most classical numerical methods, namely the finite elements method for  the Dirichlet problem, fits with our setting. For the equation $\Delta u = f,$ with Dirichlet conditions at the border, we define the set of parameters as composed by the possible functions $f$ and the set of triangulations $\mathcal{T}.$
Smoothness of the solution $u$ on $f$ is already known, but we show here that the sequence $(u_n)$ of approximations of $u$ through the finite element method is smooth for the Cauchy diffeology, with respect to the chosen triangulation and the function $f.$ For this purpose, the adapted differential geometry of the space of {{}triangulations} is described in terms of diffeologies.  
    
\section{Preliminaries on diffeology}
This section provides background on diffeology and related topics necessary for the rest of this paper.  The main reference is \cite{Igdiff}, and the reader should consult this for proofs. A complementary non exhaustive bibliography is \cite{BIgKWa2014,BN2005,BT2014,CN,CSW2014,CW2014,Don,DN2007-1,DN2007-2,FK,IgPhD,Leandre2002,Ma2006-3,Ma2013,pervova,Sou,Wa}.

\subsection{Basics of Diffeology}\label{ss:diffeology}

In this subsection we review the basics of the theory of diffeological spaces; in particular, their definition, categorical properties, as well as their induced topology. {{} The main idea of diffeologies (and Fr\"olicher spaces defined shortly after) is to replace the atlas of a classical manifold by other intrinsic objects that enable to define smoothness of mappings in a safe way, considering manifolds as a restricted class of examples. Many such settings have been developped independently. We choose these two settings because they carry nice properties such as cartesian closedness, carrying the necessary fundamental properties of e.g. calculus of variations, and also because they are very easy to use in a differential geometric way of thinking. The fundamental idea of these two settings consists in defining families of smooth maps, with mild conditions on them that ensure technical features of interest. }

\begin{Definition}[Diffeology] \label{d:diffeology}
	Let $X$ be a set.  A \textbf{parametrisation} of $X$ is a
	map of sets
	$p \colon U \to X$ where $U$ is an open subset of Euclidean space (no fixed dimension).  A \textbf{diffeology} $\p$ on $X$ is a set of
	parametrisations satisfying the following three conditions:
	\begin{enumerate}
		\item (Covering) For every $x\in X$ and every non-negative integer
		$n$, the constant function $p\colon \R^n\to\{x\}\subseteq X$ is in
		$\p$.
		\item (Locality) Let $p\colon U\to X$ be a parametrisation such that for
		every $u\in U$ there exists an open neighbourhood $V\subseteq U$ of $u$
		satisfying $p|_V\in\p$. Then $p\in\p$.
		\item (Smooth Compatibility) Let $(p\colon U\to X)\in\p$.
		Then for every $n$, every open subset $V\subseteq\R^n$, and every
		smooth map $F\colon V\to U$, we have $p\circ F\in\p$.
	\end{enumerate}
	A set $X$ equipped with a diffeology $\p$ is called a
	\textbf{diffeological space}, and is {{} denote}d by $(X,\p)$.
	When the diffeology is understood, we will drop the symbol $\p$.
	The parametrisations $p\in\p$ are called \textbf{plots}.
\end{Definition}

{{}
\noindent\textbf{Notation.} We recall that $\N^* = \{n \in \N \, | \, n \neq 0\}$ and that $\forall m \in \N^*, \N_m = \{1,...,m\} \subset \N.$}

\begin{Definition}[Diffeologically Smooth Map]\label{d:diffeolmap}
	Let $(X,\p_X)$ and $(Y,\p_Y)$ be two diffeological
	spaces, and let $F \colon X \to Y$ be a map.  Then we say that $F$ is
	\textbf{diffeologically smooth} if for any plot $p \in \p_X$,
	$$F \circ p \in \p_Y.$$
\end{Definition}

Diffeological spaces with diffeologically smooth maps form a category. This category is complete and co-complete, and forms a quasi-topos (see \cite{BH}).

\begin{Proposition}
	\cite{Sou,Igdiff} Let $(X',\p)$ be a diffeological space,
	and let $X$ be a set. Let $f:X\rightarrow X'$ be a map.
	We define $f^*(\p)$ the \textbf{pull-back diffeology} as {{} $$f^*( \p)= \left\{ p: D(p) \rightarrow X \, |   f \circ p \in \p \right\}. $$ }
\end{Proposition}

\begin{Proposition} \cite{Sou,Igdiff} Let $(X,\p)$ be a diffeological space,
	and let $X'$ be a set. Let $f:X\rightarrow X'$ be a map.
	We define $f_*(\p)$ the \textbf{push-forward diffeology} as the coarsest (i.e. the smallest for inclusion) among the diffologies
	on $X'$, which contains $f \circ \p.$ \end{Proposition}  

\begin{Definition}
Let $(X,\p)$ and $(X',\p')$
be two diffeological spaces.  A map
$f : X \rightarrow X'$
is  called  a
subduction
if  $\p' = f_*(\p).$ \end{Definition}
.
In particular, we have the following constructions.

\begin{Definition}[Product Diffeology]\label{d:diffeol product}
	Let $\{(X_i,\p_i)\}_{i\in I}$ be a family of diffeological spaces.  Then the \textbf{product diffeology} $\p$ on $X=\prod_{i\in I}X_i$ contains a parametrisation $p\colon U\to X$ as a plot if for every $i\in I$, the map $\pi_i\circ p$ is in $\p_i$.  Here, $\pi_i$ is the canonical projection map $X\to X_i$. 
\end{Definition}

In other words, in last definition, $\p = \cap_{i \in I} \pi_i^*(\p_i)$ and each $\pi_i$ is a subduction. 

\begin{Definition}[Subset Diffeology]\label{d:diffeol subset}
	Let $(X,\p)$ be a diffeological space, and let $Y\subseteq X$.  Then $Y$ comes equipped with the \textbf{subset diffeology}, which is the set of all plots in $\p$ with image in $Y$.
\end{Definition}
If $X$ is a smooth manifolds, finite or infinite dimensional, modelled on a complete locally convex topological vector space, we define the \textbf{nebulae diffeology}
$$\p(X) = \left\{ p \in C^\infty(O,X) \hbox{ (in the usual sense) }| O \hbox{ is open in } \R^d, d \in \N^* \right\}.$$
\subsection{Fr\"olicher spaces}
\begin{Definition} $\bullet$ A \textbf{Fr\"olicher} space is a triple
	$(X,\F,\mcc)$ such that
	
	- $\mcc$ is a set of paths $\R\rightarrow X$,
	
	- A function $f:X\rightarrow\R$ is in $\F$ if and only if for any
	$c\in\mcc$, $f\circ c\in C^{\infty}(\R,\R)$;
	
	- A path $c:\R\rightarrow X$ is in $\mcc$ (i.e. is a \textbf{contour})
	if and only if for any $f\in\F$, $f\circ c\in C^{\infty}(\R,\R)$.
	
	\vskip 5pt $\bullet$ Let $(X,\F,\mcc)$ et $(X',\F',\mcc')$ be two
	Fr\"olicher spaces, a map $f:X\rightarrow X'$ is \textbf{differentiable}
	(=smooth) if and only if one of the following equivalent conditions is fulfilled:
	\begin{itemize}
		\item $\F'\circ f\circ\mcc\subset C^{\infty}(\R,\R)$
	\item $f \circ \C \subset \C'$
	\item $\F'\circ f \subset  \F$ 
	\end{itemize}
\end{Definition}

Any family of maps $\F_{g}$ from $X$ to $\R$ generate a Fr\"olicher
structure $(X,\F,\mcc)$, setting \cite{KM}:

- $\mcc=\{c:\R\rightarrow X\hbox{ such that }\F_{g}\circ c\subset C^{\infty}(\R,\R)\}$

- $\F=\{f:X\rightarrow\R\hbox{ such that }f\circ\mcc\subset C^{\infty}(\R,\R)\}.$

One easily see that $\F_{g}\subset\F$. This notion will be useful
in the sequel to describe in a simple way a Fr\"olicher structure.
A Fr\"olicher space carries a natural topology,
which is the pull-back topology of $\R$ via $\F$. In the case of
a finite dimensional differentiable manifold, the underlying topology
of the Fr\"olicher structure is the same as the manifold topology. In
the infinite dimensional case, these two topologies differ very often.

 Let us now compare Fr\"olicher spaces with diffeological spaces, with the following diffeology {{}$\p_\infty(\F)$} called "nebulae":
{{}
{Let }$O${ be an open subset of a Euclidian space; } $$\p_\infty(\F)_O=
\coprod_{p\in\N}\{\, f : O \rightarrow X; \, \F \circ f \subset C^\infty(O,\R) \quad \hbox{(in
	the usual sense)}\}$$
and 
$$ \p_\infty(\F) = \bigcup_O \p_\infty(\F)_O,$$
where the latter union is extended over all open sets $O \subset \R^n$ for $n \in \N^*.$ 
}
With this construction, we get a natural diffeology when
$X$ is a Fr\"olicher space. In this case, one can easily show the following:
\begin{Proposition} \label{Frodiff} \cite{Ma2006-3} 
	Let$(X,\F,\mcc)$
	and $(X',\F',\mcc')$ be two Fr\"olicher spaces. A map $f:X\rightarrow X'$
	is smooth in the sense of Fr\"olicher if and only if it is smooth for
	the underlying nebulae diffeologies. \end{Proposition}

Thus, we can also state intuitively:
\vskip 12pt
\begin{tabular}{ccccc}
	smooth manifold  & $\Rightarrow$  & Fr\"olicher space  & $\Rightarrow$  & Diffeological space\tabularnewline
\end{tabular}
\vskip 12pt
With this construction, any complete locally convex topological vector space is a diffeological vector space, that is, a vector space for which addition and scalar multiplication is smooth. The same way, any finite or infinite dimensional manifold $X$ has a nebulae diffeology, which fully determines smooth functions from or with values in $X.$We now finish the comparison of the notions of diffeological and Fr\"olicher 
space following mostly \cite{Ma2006-3,Wa}:

\begin{Theorem} \label{compl-fro}
	Let $(X,\p)$ be a diffeological space. There exists a unique Fr\"olicher structure
	$(X, \F_\p, \mcc_\p)$ on $X$ such that for any Fr\"olicher structure $(X,\F,\mcc)$ on $X,$ these two equivalent conditions are fulfilled:
	
	(i)  the canonical inclusion is smooth in the sense of Fr\"olicher $(X, \F_\p, \mcc_\p) \rightarrow (X, \F, \mcc)$
	
	(ii) the canonical inclusion is smooth in the sense of diffeologies $(X,\p) \rightarrow (X, \p_\infty(\F)).$ 
	
	\noindent Moreover, $\F_\p$ is generated by the family 
	$$\F_0=\lbrace f : X \rightarrow \R \hbox{ smooth for the 
		usual diffeology of } \R \rbrace.$$
	{{} We call \textbf{Fr\"olicher completion} of $\p$ the Fr"olicher structure $(X, \F_\p, \mcc_\p).$}
\end{Theorem}

\noindent
\textbf{Proof.}
Let $(X,\F,\mcc)$ be a Fr\"olicher structure satisfying \textit{(ii)}. 
Let $p\in P$ of domain $O$. $\F \circ p \in C^\infty(O,\R)$ in the usual sense.  
Hence, if $(X,\F_\p, \mcc_\p)$is the Fr\"olicher structure on $X$ generated by the 
set of smooth maps $(X,\p)\rightarrow \R,$ we have two smooth inclusions 
$$ (X,\p) \rightarrow (X,\p_\infty(\F_\p)) \hbox{ in the sense of diffeologies }$$
and
$$ (X, \F_\p, \mcc_\p) \rightarrow (X,\F,\mcc) \hbox{ in the sense of Fr\"olicher. }$$
Proposition \ref{Frodiff} ends the proof. \quad \qed

\begin{Definition} \cite{Wa}
	A \textbf{reflexive} diffeological space is a diffeological space $(X,\p)$ such that $\p = \p_\infty(\F_\p).$
\end{Definition}

\begin{Theorem} \cite{Wa}
	The category of Fr\"olicher spaces is exactly the category of reflexive diffeological spaces.
\end{Theorem}

This last theorem allows us to make no difference between Fr\"olicher spaces and reflexive diffeological spaces. 
We shall call them Fr\"olicher spaces, even when working with their underlying diffeologies.
\vskip 12pt A deeper analysis of these implications has been given
in \cite{Wa}. The next remark is inspired on this work and on
\cite{Ma2006-3}; it is based on \cite[p.26, Boman's theorem]{KM}.
\begin{rem}
	We notice that the set of contours $\C$ of the Fr\"olicher space
	$(X,\F,\C)$ \textbf{does not} give us a diffeology, because a diffelogy
	needs to be stable under restriction of domains. In the case of paths in
	$\C$ the domain is always $\R.$ However, $\C$ defines a ``minimal diffeology''
	$\p_1(\F)$ whose plots are smooth parameterizations which are locally of the
	type $c \circ g,$ where $g \in \p_\infty(\R)$  and $c \in \C.$ Within this setting,
	a  map $f : (X,\F,\C) \rightarrow (X',\F',\C')$ is smooth if and only if it is smooth  $(X,\p_\infty(\F)) \rightarrow (X',\p_\infty(\F')) $ or equivalently smooth  .$(X,\p_1(\F)) \rightarrow (X',\p_1(\F')) $ 
\end{rem}
We apply the results on product diffeologies to the case of Fr\"olicher spaces and we derive very easily, (compare with e.g. \cite{KM}) the following:

\begin{Proposition} \label{prod2} Let $(X,\F,\C)$
	and $(X',\F',\C')$ be two Fr\"olicher spaces equipped with their natural
	diffeologies $\p$ and $\p'$ . There is a natural structure of Fr\"olicher space
	on $X\times X'$ which contours $\C\times\C'$ are the 1-plots of
	$\p\times\p'$. \end{Proposition}

We can even state the result above for the case of infinite products;
we simply take cartesian products of the plots or of the contours.
We also remark that given an algebraic structure, we can define a
corresponding compatible diffeological structure. For example, a
$\R-$vector space equipped with a diffeology is called a
diffeological vector space if addition and scalar multiplication
are smooth (with respect to the canonical diffeology on $\R$), see \cite{Igdiff,pervova,pervova2}. An
analogous definition holds for Fr\"olicher vector spaces. Other
examples will arise in the rest of the text.

\begin{rem} \label{comp}
	Fr\"olicher, $c^\infty$ and Gateaux smoothness are the same notion
	if we restrict to a Fr\'echet context, see \cite[Theorem 4.11]{KM}.
	Indeed, for a smooth map $f : (F, \p_1(F)) \rightarrow \R$ defined
	on a Fr\'echet space with its 1-dimensional diffeology, we have
	that $\forall (x,h) \in F^2,$ the map $t \mapsto f(x + th)$ is
	smooth as a classical map in $\C^\infty(\R,\R).$ And hence, it is
	Gateaux smooth. The converse is obvious.
\end{rem}

\subsection{Quotient and subsets}

We give here only the results that will be used in the sequel.

We have now the tools needed to describe the diffeology on a quotient:

\begin{Proposition} \label{quotient} let $(X,\p)$ b a diffeological
	space and $\rel$ an equivalence relation on $X$. Then, there is
	a natural diffeology on $X/\rel$, {{} denote}d by $\p/\rel$, defined as
	the push-forward diffeology on $X/\rel$ by the quotient projection
	$X\rightarrow X/\rel$. \end{Proposition}

Given a subset $X_{0}\subset X$, where $X$ is a Fr\"olicher space
or a diffeological space, we can define on subset structure on $X_{0}$,
induced by $X$.

$\bullet$ If $X$ is equipped with a diffeology $\p$, we can define
a diffeology $\p_{0}$ on $X_{0},$ called \textbf{subset
	diffeology} \cite{Sou,Igdiff} setting \[ \p_{0}=\lbrace p\in\p
\hbox{ such that the image of }p\hbox{ is a subset of
}X_{0}\rbrace.\]

\begin{example}
    Let $X$ be a diffeological space. Let us {{} denote} by $$X^\infty = \left\{ (x_n)_{n \in \N} \in X^\N \, | \, \{n \, | x_n \neq 0 \} \hbox{ is a finite set}\right\} $$ Then this is a diffeological space, as a subset of $X^\N.$
    \end{example}

$\bullet$ If $(X,\F,\C)$ is a Fr\"olicher space, we take as a generating
set of maps $\F_{g}$ on $X_{0}$ the restrictions of the maps $f\in\F$.
In that case, the contours (resp. the induced diffeology) on $X_{0}$
are the contours (resp. the plots) on $X$ which image is a subset
of $X_{0}$.

\subsection{Projective limits and vector pseudo-bundles}

Let us now give the description of what happens
for projective limits of Fr\"olicher and diffeological spaces.

\begin{Proposition} \label{froproj} Let $\Lambda$ be an infinite
	set of indexes which can even be uncountable.
	
	$\bullet$ Let $\lbrace(X_{\alpha},\p_{\alpha})\rbrace_{\alpha\in\Lambda}$
	be a family of diffeological spaces indexed by $\Lambda$ totally
	ordered for inclusion,  with $(i_{\beta,\alpha}: X_\alpha \rightarrow X_\beta)_{(\alpha,\beta)\in\Lambda^{2}}$
	the family of inclusion maps which are assumed smooth maps. If $X=\bigcap_{\alpha\in\Lambda}X_{\alpha},$ then
	$X$ carries the \textbf{projective diffeology} $\p$ which is the
	pull-back of the diffeologies $\p_{\alpha}$ of each $X_{\alpha}$
	via the family of inclusion maps $(f_{\alpha}:X \rightarrow X_\alpha)_{\alpha\in\Lambda}.$ The diffeology
	$\p$ is made of plots $g:O\rightarrow X$ such that for each $\alpha\in\Lambda,$
	\[
	f_{\alpha}\circ g\in\p_{\alpha}.\]
	This is the biggest diffeology for which the maps $f_{\alpha}$ are
	smooth.
	
	$\bullet$ We have the same kind of property for Fr\"olicher spaces: let $\lbrace(X_{\alpha},\F_{\alpha},\C_{\alpha})\rbrace_{\alpha\in\Lambda}$
	be a family of Fr\"olicher spaces indexed by $\Lambda$, a non-empty set
	totally ordered
	for inclusion.  With the sae notations,there is a natural structure of Fr\"olicher space on  $X=\bigcap_{\alpha\in\Lambda}X_{\alpha},$ for
	which the corresponding contours
	\[  \C=\bigcap_{\alpha\in\Lambda}\C_{\alpha}   \]
	are some 1-plots of $\p=\bigcap_{\alpha\in\Lambda}\p_{\alpha}.$ A
	generating set of functions for this Fr\"olicher space is the set of
	maps of the type:
	$$
	\bigcup_{\alpha \in \Lambda} \F_{\alpha} \circ f_{\alpha}.$$
\end{Proposition}
{{} Let us now have a precise look at the notion of fiber bundle in clasiacl (finite dimensional) fiber bundles. Fiber bundles, in the context of smooth finite dimensional manifolds, are defined by \begin{itemize}
		\item a smooth manifold  $E$ called total space
		\item a smooth manifold  $X$ called base space
		\item a smooth submersion $\pi: E \rightarrow X$ called fiber bundle projection
		\item a smooth manifold $F$ called typical fiber, because $\forall x \in X,  \pi^{-1}(x)$ is a smooth submanifold of $E$ diffeomorphic to $F.$
		\item a smooth atlas on $X,$ with domains $U \subset X$ such that $\pi^{-1}(U)$ is an open submanifold of $E$ diffeomorphic to $U \times F.$ We the get a system of local trivializations of the fiber bundle.
\end{itemize}
By the way, in order to be complete, a smooth fiber bundle should be the quadruple data $(E,X,F,\pi)$ (because the definition of $\pi$ and of $X$ enables to find systems of local trivializations). For short, this quadruple setting is often {{} denote}d by the projection map $\pi: E\rightarrow X.$

There exists some diffeological spaces which carry no atlas, so that, the condition of having a system of smooth trivializations in a generalization of the notion of fiber bundles is not a priori necessary, even if this condition, which is additional, enables interesting technical aspects \cite[pages 194-195]{MW2017}. So that, in a general setting, we do not need to assume the existence of local trivializations.}
Now, following \cite{pervova}, in which the ideas from \cite[last section]{Sou} have been devoloped to vector spaces, the notion of quantum structure has been introduced in \cite{Sou} as a generalization of principal bundles, and the notion of vector pseudo-bundle in \cite{pervova}.The common idea consist in the description of fibered objects made of a total (diffeological) space $E,$ over a diffeological space $X$ and with a canonical smooth {{} bundle} projection $\pi: E \rightarrow X$ such as, $\forall x \in X,$ $\pi^{-1}(x)$ is endowed with a (smooth) algebraic structure, but for which we do not assume the existence of a system of local trivialization. 
\begin{enumerate}
	\item For a diffeological vector pseudo-bundle, the fibers $\pi^{-1}(x)$ are assumed diffeological vector spaces, i.e. vector spaces where addition and multiplication over a diffeological field of scalars (e.g. $\R$ or $\mathbb{C}$) is smooth. We notice that \cite{pervova} only deals with finite dimensional vector spaces.
	\item For a so-called ``structure quantique'' (i.e. ``quantum structure'') following the terminology of \cite{Sou}, a diffeological group $G$ is acting on the right, smoothly and freely on a diffeological space $E$. The space of orbits $X=E/G$ defines the base of the quantum structure $\pi: E \rightarrow X,$ which generalize the notion of principal bundle by not assuming the existence of local trivialization. In this picture, each fiber $\pi^{-1}(x)$ is isomorphic to $G.$
\end{enumerate}
From these two examples, we can generalize the picture. 
\begin{Definition}\label{pseu-fib}
	Let $E$ and $X$ be two diffeological spaces and let $\pi:E\rightarrow X$ be a smooth surjective map. Then $(E,\pi,X)$ is a \textbf{diffeological fiber pseudo-bundle} if and only if  $\pi$ is a subduction. 
\end{Definition}
{{} Let us precise that we do not assume that there exists a typical fiber, in coherence with Pervova's diffeological vector pseudo-bundles. We
}
 can give the following definitions:
\begin{Definition}
	Let $\pi:E\rightarrow X$ be a diffeological fiber pseudo-bundle. Then:
	\begin{enumerate}
		\item Let $\mathbb{K}$ be a diffeological field. $\pi:E\rightarrow X$ is a diffeological $\mathbb{K}-$vector pseudo-bundle if there exists: 
		\begin{itemize}
			\item a smooth fiberwise map $.\, :\mathbb{K} \times E \rightarrow E,$
			\item a smooth fiberwise map $+:E^{(2)} \rightarrow E$ where $$E^{(2)} = \coprod_{x \in X} \{(u,v) \in E^2\, | \, (u,v)\in \pi^{-1}(x)\}$$ equipped by the pull-back diffeology of the canonical map $E^{(2)} \rightarrow E^2,$
		\end{itemize}
		such that $\forall x \in X, $ $(\pi^{-1}(x),+,.)$ is a diffeological $\mathbb{K}-$vector bundle.
		\item $\pi:E\rightarrow X$ is a \textbf{diffeological gauge pseudo-bundle} if there exists \begin{itemize}
			\item a smooth fiberwise involutive map ${(.)}^{-1}\, E \rightarrow E,$
			\item a smooth fiberwise map $.\,:E^{(2)} \rightarrow E$ 
		\end{itemize}
		such that $\forall x \in X, $ $(\pi^{-1}(x),\, .\, )$ is a diffeological group with inverse map $(.)^{-1}.$
		\item  $\pi:E\rightarrow X$ is a \textbf{diffeological principal pseudo-bundle} if there exists a diffeological gauge pseudo-bundle $\pi': E' \rightarrow X$ such that, considering $$E\times_X E' = \coprod_{x \in X} \{(u,v)\, | \, (u,v)\in \pi^{-1}(x)\times \pi'^{-1}(x)\}$$ equipped by the pull-back diffeology of the canonical map $E\times_X E' \rightarrow E\times E',$ there exists a smooth map $E\times_X E' \rightarrow E$ which restricts fiberwise to a smooth free and transitive right-action $$\pi^{-1}(x) \times \pi'^{-1}(x) \rightarrow \pi^{-1}(x).$$ 
		\item $\pi:E\rightarrow X$ is a \textbf{Souriau quantum structure} if it is a diffeological principal pseudo-bundle with diffeological gauge (pseudo-)bundle $X\times G \rightarrow X.$ 
	\end{enumerate}
\end{Definition}
\section{Tangent spaces, diffeology and group of diffeomorphisms} \label{diff} \label{T2}

There are actually two main definitions, (2) and (3) below, of the tangent space of a diffeological space $(X,\p),$ while definition (1) of the tangent cone will be used in the sequel. These definitions are very similar to the defnitions of the two tangent spaces in the $c^\infty-$setting given in \cite{KM}.

\begin{enumerate}
	\item the \textbf{internal tangent cone}. This construction is extended to the category of diffeologies from the definitions in \cite{DN2007-1} on Fr\"olicher spaces, and very similar to the kinematic tangent space defined in \cite{KM}. 
	For each $x\in X,$ we consider $$C_{x}=\{c \in C^\infty(\R,X)| c(0) = x\}$$ 
	and take the equivalence relation $\mathcal{R}$ given by $$c\mathcal{R}c' \Leftrightarrow \forall 
	f \in C^\infty(X,\R), \partial_t(f \circ c)|_{t = 0} = \partial_t(f \circ c')|_{t = 0}.$$
	The internal tangent cone at $x$ is the quotient
	$$^iT_xX = C_x / \mathcal{R}.$$ If $X = \partial_tc(t)|_{t=0} \in {}^iT_X, $ we define the simplified notation  $$Df(X) = \partial_t(f \circ c)|_{t = 0}.$$
	Under these constructions, we can define the total space of internal tangent cones $$ {}^{i}TX = \coprod_{x \in X} {}^{i}T_xX$$ with canonical projection $\pi: u \in {}^{i}T_xX \mapsto x,$ and equipped with the diffeology defined by the plots $p : D(p) \subset \R^n \rightarrow {}^{i}TX$ defined through the plots $p': (t,z)\in\R \times D(p) \rightarrow X \in \p$ by $$p(z) = \partial_t p'(t,z)|_{t=0}.$$ With this construction, ${{}^{i}}TX$ is a diffeological fiber pseudo-bundle.
	This construction can be criticized because, depending on the base diffeology $\p$ on $X,$, there can exist ``too few'' maps $f \in C^\infty(X,\R)$ and hence some germs of paths $c$ cannot be separated via evaluations by smooth functions $f.$  
	\item The \textbf{internal tangent space} is defined in \cite{He1995,CW2014}, based on germs of paths. {{} This second definition is necessary, and the internal tangent space differ from the internal tangent cone. Indeed, spaces of germs do not carry intrinsically a structure of abelian group. This remark was first formulated in the context of Fr\"olicher spaces, see \cite{DN2007-1}, and see \cite{CW2014} for the generalization to diffeologies. For this reason, one can complete the tangent cone into a vector space, called internal tangent space. This was performed in \cite{CW2014} via mild considerations on colimits in categories.}
	\item the \textbf{external} tangent space $^eTX,$ defined as the set of derivations on $C^\infty(X,\R).$ \cite{KM,Igdiff}.
\end{enumerate}

\noindent
For finite dimensional {{}manifolds}, definitions (1), (2) and (3) coincide.
However, to our best knowledge, there is still one aspect that is not investigated yet among all possible ways to generalize tangent spaces. For this purpose, we need to recall the following definitions from \cite{Igdiff}: \begin{Definition}

	Let $(X,\p)$ and $(X',\p')$ be two diffeological spaces. Let $S \subset C^\infty(X,X')$ be a set of smooth maps. The \textbf{{
			functional} diffeology} on $S$ is the diffeology $\p_S$
	made of plots
	$$ \rho : D(\rho) \subset \R^k \rightarrow S$$
	such that, 
	for each $p \in \p, $
	the maps $\Phi_{\rho, p}: (x,y) \in D(p)\times D(\rho) \mapsto \rho(y)(x) \in X'$ are plots of $\p'.$
\end{Definition}

\noindent
With this definition, we have the classical fundamental property for calculus of variations and for composition:

\begin{Proposition} \label{functf}\cite{Igdiff}
	Let $X,Y,Z$ be diffeological spaces
	\begin{enumerate}
		\item
	$$C^\infty(X\times Y,Z) = C^\infty(X,C^\infty(Y,Z)) = C^\infty(Y,C^\infty(X,Z))$$as diffeological spaces equipped with {
		functional} diffeologies.
	\item The composition map 
	$$C^\infty(X,Y) \times C^\infty(Y,Z) \rightarrow C^\infty(X,Z)$$ is smooth.  
\end{enumerate}
	
\end{Proposition}
Let us now investigate tangent spaces from the viewpoint of diffeomorphisms. On a finite dimensional manifold $M,$ the Lie algebra of the ILH Lie group of diffeomorphisms \cite{Om}, defined as the tangent cone at the identity map, is the space of smooth vector fields, i.e. smooth sections of $TM.$ On a non-compact, locally compact manifold, the situation is quite similar \cite{Ma2013-2} while the group of diffeomorphisms is no longer a Fr\'echet manifold but a Fr\"olicher Lie group. On these groups, the underlying diffeology is the functional diffeology, as well as for the more general definition of $Diff(X)$ when $X$ is a diffeological space \cite{Igdiff}.

We now get the necessary material to give the following definition. In the rest of this section, $(X,\F,\C)$ is a Fr\"olicher space. 

\begin{Definition} \label{def:dtangent}We use here the notations that we  used before for the definition of the internal tangent cone. 
	Let ${}^{d}T_xX$ be the subset of ${}^iT_xX$ defined by $$ {}^dT_xX = {}^dC_x /\mathcal{R}$$ with $${}^dC_x =\left\{ c \in C_x | \exists \gamma \in C^\infty(\R,Diff(X)), c(.)=\gamma(.)(x) \hbox{ and } \gamma(0) = Id_X  \right\}$$
\end{Definition}

Through this definition, ${}^dT_xX$ is intrinsically linked with the tangent space at the identity ${}^iT_{Id_X}Diff(X)$ described in \cite{Les} for any diffeological group (i.e. group equipped with a diffeology which makes composition and inversion smooth), see e.g. \cite{MR2016}, as smooth vector space.
\begin{rem} \label{rq24}
	Let $\gamma \in C^\infty(\R,Diff(X))$ such that $\gamma(0)(x)=x.$ Then $\lambda (x) = (\gamma(0))^{-1} \circ \gamma(.)(x)$ defines a smooth path $\lambda \in {}^dC_x.$ {{} Consequently,}  $${}^dC_x =\left\{ c \in C_x | \exists \gamma \in C^\infty(\R,Diff(X)), c(.)=\gamma(.)(x) \hbox{ and } \gamma(0)= Id_X  \right\}$$ 
\end{rem}
\begin{Definition}
	Let $X$ be a Fr\"olicher space.
	we define, by $$ {}^dTX = \coprod_{x \in X} {}^dT_xX$$
	the diff-tangent bundle of $X.$ 
\end{Definition} 

By the way, we can get easily the following observations: 

\begin{Proposition} \label{prop:pdiff}
	Let $(X,\p)$ be a reflexive diffeological space, and let $\p_{Diff}$ be the functional diffeology on $Diff(X).$
	\begin{enumerate}
		\item \label{dt1} There exists a diffeology $\p(Diff) \subset \p$ {{} which is} generated by the family of push-forward diffeologies : 
		$$ \left\{ (ev_x)_{*}(\p_{Diff}) \, | \, x \in X  \right\}.$$
		\item \label{dt2} $\forall x \in X, {}^dT_xX$ is the internal tangent cone of $(X,\p(Diff))$ at $x.$
		\item \label{dt3} $\forall x \in X, {}^dT_xX$ is a diffeological vector space
		\item \label{dt4} The total diff-tangent space $${}^{d}TX = \coprod_{x \in X} {}^{d}TX \subset {}^{i}TX $$ is a vector pseudo-bundle for the subset diffeology inherited from ${}^{i}TX $ and also for the diffeology of internal tangent space of $(X,\p_{Diff}).$
	\end{enumerate}
\end{Proposition}
\begin{proof}
(\ref{dt1}) is a consequence of the definition of push-forward diffeologies the following way: the family $$\{ \p \hbox{ diffeology on } X \, | \, \forall x \in X, \, (ev_x)_{*}(\p_{Diff}) \subset \p\}$$ has a minimal element by Zorn Lemma. 

(\ref{dt2}) follows from remark \ref{rq24}.

(\ref{dt3}): The diffeology $\p(Diff)$ coincides with the diffeology made of plots which are locally of the form $ev_x \circ p,$ where $x \in X$ and $p$ is a plot of the diffeology of $Diff(X).$
We have that $^{i}T_{Id}Diff(X)$ is a diffeological vector space, following \cite{Les}. This relation follows from the differentiation of the multiplication of the group: given two paths $\gamma_1, \gamma_2$ in $C^\infty(\R,Diff(X)),$ with $\gamma_1(0)=\gamma_2(0)=Id,$ if $X_i = \partial_t\gamma_i(0)$ for $i \in \{1{{},}2\},$ then $$ X_1 + X_2 = \partial_t (\gamma_1 . \gamma_2) (0).$$ Reading locally plots in $\p(Diff),$ we can consider only plots of the for $ev_x \circ p,$ where $p$ is a plot in $Diff(X)$ such that $p(0) = Id_X.$ By the way the vector space structure on $^{d}T_xX$ is inherited from  $^{i}T_{Id}Diff(X)$ via evaluation maps.

In order to finish to check (\ref{dt3}), we prove directly $(\ref{dt4})$ by describing its diffeology.

For this, we consider $$C^\infty_0(\R, Diff(X)) = \left\{ \gamma \in C^\infty_0(\R, Diff(X)) \, | \, \gamma(0)=Id_X \right\}.$$
Let $^{d}C = \coprod_{x \in X} ^{d}C_x.$ The total evaluation map \begin{eqnarray*}
ev : & X \times C^\infty_0(\R, Diff(X)) \rightarrow & ^{d}C \\
& (x,\gamma) \mapsto & ev_x \circ \gamma
\end{eqnarray*}
is fiberwise (over $X$), and onto. By the way we get a diffeology on $ ^{d}C$ which is the push-forward diffeology of $X \times C^\infty_0(\R, Diff(X))$ by $ev.$ Passing to the quotient, we get a diffeology on $^{d}TX$ which makes each fiber $^{d}T_xX$ a diffeological vector space trivially.

\end{proof}
\begin{example}
Let us consider $X \subset \R^2$ defined by $$ X =\left\{ (x,y) \in \R^2 \, | \, xy=0  \right\}.$$
$\R^2$ is a Fr\"olicher space (equipped with its nebulae diffeology $\p_\infty(\R^2)$) and $X$ has a subset diffeology made of plots of three types: 

\begin{itemize}
\item plots of the subset diffeology of $X_1 = \{ (x,y) \in \R^2 \, | \, y = 0 \}$ which is (diffeologically) isomorphic to $(\R, \p_\infty(\R))$
\item plots of the subset diffeology of $X_2 = \{ (x,y) \in \R^2 \, | \, x = 0 \}$ which is (diffeologically) isomorphic to $(\R, \p_\infty(\R))$
\item plots which are locally in $\p_\infty(X_1)$ or $\p_\infty(X_2),$ obtained by gluing along $X_1 \cap X_2 = \{(0{{},}0)\}$ where plots have to be stationary. 
\end{itemize}
Let $$\F_1 = \left\{ f \in C^\infty(\R^2,\R) \, | \, \forall (x,y,y') \in \R^3, \, f(x,y)=f(x,y') \right\}$$
and let $$\F_2 = \left\{ f \in C^\infty(\R^2,\R) \, | \, \forall (x,y,x') \in \R^3, \, f(x,y)=f(x',y) \right\}$$
The subset diffeology of $X$ is generated by $\F_1 \cup \F_2,$ i.e. p is a plot of this fiffeology if and only if $ (\F_1 \cup \F_2) \circ p \subset C^\infty(D(p),\R).$ So that it is reflexive.  

Let us now highlight the internal tangent cone.

\begin{itemize}
\item $\forall x \in \R^*,$ $^{i}T_{(x,0) }X$ and $^{i}T_{(0,x) }X$ are both diffeologically isomorphic to $\R$
\item The internal tangent cone at the origin $^{i}T_{(0{{},}0)} X \sim \R \cup  \R (\subset \R^2 = T_{(0{{},}0)}\R^2)$ is a cone, and which completes to $\R^2$ along the lines of the work by Christensen and Wu \cite{CW2014}.
\end{itemize}
Let us now consider $g \in Diff(X).$ $G$ is continuous for the $D-$topology of $X.$ Let us consider $z = g(0{{},}0).$ The set $X - \{(0{{},}0)\}$ has four connected components for the $D-$toplology, so that $g(X - \{(0{{},}0)\})= g(X) - \{z\}$ has also four connected components. By the way, $z = (0{{},}0)$ and hence $(0{{},}0)$ is a fixed point for any diffeomorphism $g$ of $X.$ This shows that $$^{d}T_{(0{{},}0)}X = \{0\} \neq                   {}^{i}T_{(0{{},}0)}X$$ and hence $$ ^{d}TX \neq {}^{i}TX$$ while $ ^{d}TX \subset {}^{i}TX.$
\end{example}
\section{Cauchy diffeology and smooth numerical schemes}
\subsection{The Cauchy diffeology} \label{2}

Let $X \subset Y,$ let $j: X \rightarrow Y$ the canonical inclusion. In what follows, we assume that:

\begin{itemize}
	\item $X$ is a diffeological space with diffeology $\p.$
	\item $Y$ is a $(T_2)$ sequentially complete uniform space \cite[Topologie g\'en\'erale, Chapter II]{Bou}. 
	\item $\forall p \in \p, $ $j \circ p$ is a continuous map. 
	\item $Y$ is equipped with a diffeology $\p_0.$
\end{itemize}
\vskip 6pt
The diffeology $\p_0$ is made of plots which can be not continuous for the uniform space topology. This is why we consider the following diffeology:

\vskip 6pt
- Let $\F_0 = \{ f \in C^0(Y, \R) |\forall p \in \p_0, f \circ p \hbox{ is smooth} \}$

- Let $\mcc = \mcc(\F_0)$ and $\F = \F(\mcc).$

- Let $\p' = \{p \in \p_\infty(\F)| p \hbox{ is continuous}\}.$
\vskip 6pt
Notice that if $C^0(Y, \R)$ determines the topology of $Y$ by pull-back of the topology of $\R,$ the last point is not necessary, and $\p' = \p_{{}\infty}(\F).$ 

\begin{Definition}
We {{} denote} by $\mcc(X,Y)$ the subspace of $X^\N$ made of sequence $(x_n)_{n \in \N}$ such that $\left(j (x_n)\right)_{n \in \N}$ is Cauchy in $Y.$
\end{Definition}

\noindent
When $(X,\p)=(Y,\p_0)$ as diffeological spaces, we shall use the notation $\mcc(X)$ instead of $\mcc(X,X).$
There exists two naive "natural" diffeologies on $\mcc(X):$
\begin{itemize}
	\item  As an infinite product, $X^\N$ carries a natural diffeology $\p_\N,$ and by the inclusion $\mcc(X) \subset X^\N,$ the set $\mcc(X)$ can be endowed with the subset diffeology.
	\item If $X$ is complete, $\mcc(X)$ can be equipped with the pull-back diffeology $$ \left(\lim_{n \rightarrow +\infty }\right)^*(\p').$$
\end{itemize}
 {{}The two diffeologies cannot be compared}, as it is shown in next proposition: 
 
 \begin{Proposition}
 	If there exists two distinct points $x$ and $y$ in $X$ which are connected by a smooth path $\gamma : [0;1] \rightarrow X, $ then $$\lim_{n \rightarrow +\infty}: (\mcc(X),\p_\N)\rightarrow (X,\p)$$ is not smooth.
 \end{Proposition}   
 \begin{proof}
 Let \textbf{x} be the constant Cauchy sequence which converges to $x,$ and let \textbf{y} be the constant Cauchy sequence which converges to $y.$ Let us define a path $\Gamma: [0,1]\rightarrow \mcc(X)$ by the relations, for $n\geq 0:$
 \begin{eqnarray}
 \Gamma(t)_n & = & \left\{\begin{array}{cl}
 x & \hbox{ if } t < 1-\frac{1}{n+1} \\
 \gamma\left( (n+2)(n+1)\left(t-1+\frac{1}{n+1}\right)\right) & \hbox{ if } t \in \left[ 1 - \frac{1}{n+1};1 - \frac{1}{n+2} \right] \\
 y & \hbox{ if } t> 1 - \frac{1}{n+2} \end{array}
  \right.
 \end{eqnarray}
 This path is smooth for $\p_\N,$ and 
 $$\left\{
 \begin{array}{cccl}
 &\lim_{n\rightarrow +\infty}\Gamma_n(0) & = &x \\
 \forall t\in \left]0;1\right], &\lim_{n\rightarrow +\infty}\Gamma_n(t) & = &y \end{array}
 \right.$$
 On this path, {{} $$\lim_{t \rightarrow 0^+} \left( \lim_{n\rightarrow +\infty}\Gamma_n(t) \right) \neq \lim_{n\rightarrow +\infty}\Gamma_n(0)$$
 	 which shows that the map $\lim_{n\rightarrow +\infty}$ 
 	 is not continuous for the $D-$topologies of $\p_\N$ and $\p;$ thus} the map 
 	 $\lim_{n \rightarrow +\infty}$ 
 	 is not smooth. \end{proof}
 
 \noindent
 We illustrate also that the two diffeologies are very different by the following example, already given in \cite{ERMR2017}:
 
 \begin{example} \label{exQ}
 	Let us equip $\mathbb{Q}\subset \mathbb{R}$ with the subset diffeology of $\R.$ This is the discrete diffeology. Thus the diffeology $\p_\N$ on $\mcc(\mathbb{Q})$ is the discrete diffeology, and hence  the map $\lim_{n \rightarrow +\infty}: \mcc(\mathbb{Q})\rightarrow \R$ is smooth. However, $\p_\N \neq (\lim_{n \rightarrow +\infty})^* \p_\infty(\R).$ 
 \end{example}
 
 \noindent
 Since we have motivated before  the use of the two diffeologies, and since both diffeologies cannot be compared in an efficient way, the following definition becomes natural:
 
 \begin{Definition}
 	Let $X$ and $Y$ as above. We define the \textbf{Cauchy diffeology} on $\C(X,Y),$ {{} denote}d $\p_\mcc,$ as $$\p_\mcc = \p_\N \cap \left(\lim_{n \rightarrow +\infty }\right)^*(\p),$$ where $\lim$ is the $Y-$limit and $\p_\N$ is the subset diffeology in $\C(X,Y)$ inherited from $X^\N.$
 \end{Definition}

\noindent
From example \ref{exQ} we see that the ambient space $X$ needs to have ``enough'' plots to make the Cauchy diffeology interesting. There exist already a non-trivial example of Cauchy diffeology in the literature, in the setting of an ultrametric completion: 

\begin{example} \textbf{Diffeologies compatible with a valuation and ultrametric completion.}
	\cite{ERMR2017} Let $(X,val)$ be a $\mathbb{K}-$algebra with valuation. assuming that $\mathbb{K}$ is equipped with the discrete diffeology, let $\p$ be a diffeology on $X$ such that  
	\begin{enumerate}
		\item addition and multiplication $A \times A \rightarrow A$ are smooth,
		\item scalar multiplication $\mathbb{K} \times A \rightarrow A$ is smooth
		\item Let $p \in \Z$ and let $X_p = \{ x \in X \, |\, val(x)\geq p\}.$ Let us equip $X/X_p$ by the quotient diffeology and let us equip the ultrametric completion $\hat{X}$ with the pull-back diffeology of the projection maps $\pi_p: \hat{X} \rightarrow X / X_p$ for $p \in \Z.$ Then the diffeology on $X$ is the pull-back of diffeology on $\hat{X}$ through the canonical inclusion $X \rightarrow \hat{X}.$  
	\end{enumerate}
Then, the map $\lim_{n \rightarrow +\infty}: \C(X) \rightarrow \hat{X}$ is smooth \cite[Theorem 2.l]{ERMR2017}, in other words, the Cauchy diffeology coincides with the subset diffeology of $\C(X) \subset X^\N,$ i.e. $\p_\C = \p_\N.$
\end{example}

\noindent
Let us give another example, where the Cauchy diffeology is  more complex.

\begin{example} \textbf{Cauchy diffeology on $\mathbb{P}(X).$}
	Let $X$ be a compact metric space, equipped with its Borel $\sigma-$algebra. The space of probabilities on $X,$ $\mathbb{P}(X)$ is a convex subspace of the space of measures on $X$  and the set of extremals (the ``border'') of $\mathbb{P}(U)$
	is the set of Dirac measures $\delta(X)=\{\delta_x; x \in X\}.$  The Monte Carlo method is based on isobarycentres of Dirac measures, and any probability measure $\mu \in \mathbb{P}(X)$ is a limit of isobarycentres of Dirac measures for the vague topology, in other words,
	$$\forall \mu \in \mathbb{P}(X), \exists (x_n)\in X^\N, \forall f \in C^0(X,\R), \int_X f d\mu = \lim_{n \rightarrow +\infty}\frac{1}{n+1}\sum_{k=0}^n f(x_k).$$
	Conversely, each sequence of Dirac measures generate this way a probability measure.
	Equipped with the Prokhorov metric, $\mathbb{P}(X)$ is a compact metric space \cite{Prok}, {{} see e.g. \cite{Bil1999,Dud1968}}. If $X$ is a diffeological space, the family $C^\infty(X,\R)$ generates a structure of Fr\"olicher space, and hence we can define the diffeology $\p' $ on $\mathbb{P}(X)$ following the notations of the section. Thus, we get a Cauchy diffeology on $\mcc(\delta(X),\mathbb{P}(X)).$ {{}
		This diffeology contains: 
		\begin{enumerate}
			\item smooth paths for the diffeology $\p_\N$ $$\gamma : [0;1]\rightarrow (x_0(t),...,x_n(t),...)$$
			with values in sequences that are stationnary after a fixed rank , i.e. such that $\exists N \in \N, \forall n > N, \forall t,  x_n(t) = x_N(t).$ Thus the Cauchy diffeology is not the discrete diffeology.
			\item paths that are not targeted in sequences stationnary after finite rank, such as, when $U=[0;1],$ 
			$t \in [0;1] \mapsto (x_n(t))$ with $$x_n(t) = \frac{t}{n+1}.$$   
		\end{enumerate} 
	So that, this diffeology contains a rich class of smooth paths.}  
\end{example}

\noindent
This example gives a first motivation for the following theorem:

\begin{Proposition}
	Let $(X,\p_X)$ and $(Y,\p_Y)$ be diffeological vector spaces, such that $X \subset Y$ with smooth inclusion, and such that $Y$ is equipped with a translation-invariant metric which generates a topology $\tau$ which is weaker than its $D-$topology. Assume also that $X$ is dense in $Y$ for $\tau.$ Then $(\mcc(X), \p_\mcc)$ is a diffeological vector space.   
\end{Proposition}

\begin{proof}Since $\mcc(X)$ is a vector space, it is a diffeological vector space for $\p_\N.$ Since $\lim_{n \rightarrow +\infty}$ is a linear map, $\left( \lim_{n \rightarrow +\infty}\right)^*(\p_Y)$ is also a diffeology of diffeological vector space on $\mcc(X)$. Thus, so is $\p_\mcc.$ \end{proof}
 
 \subsection{A theory of smooth numerical schemes} \label{3}
 We omit in this section to precise the diffeologies under consideration when these are obvious ones: nebulae (reflexive) diffeologies for {{} locally convex topological vector spaces (LCTVS) or Fr\'echet spaces or manifolds, and arbitrary,} fixed diffeologies for diffeological spaces.

 Let {{} $X$ and $Z$ be a (LCTVS) and let $Y$ be Fr\'echet spaces. We insist here on the fact that only completeness of $Y$ is necessary for the Cauchy diffeology on $\mcc(X,Y).$} Let $Q$ be a diffeological space of parameters. 
 
 \begin{itemize}
 	\item {{} Assume that the inclusion map $X \rightarrow Y$ is smooth.}
 	\item Let us consider the space of Cauchy sequences $\mcc(X,Y)$ that are Cauchy sequences in $X$ with respect to the uniform structure on $Y.$
 \end{itemize}
 
 \begin{Definition}

 A \textbf{smooth functional equation} is defined by a smooth map $F : X\times Q \rightarrow Z$ and by the condition
 \begin{equation} \label{eq}
 F(u,q)=0
 \end{equation}
 {{} The set $Num_F(Y)$ of $Y-$\textbf{smooth numerical schemes} is the set of smooth maps $$x: Q \rightarrow \mcc(X,Y)$$
 such that, if $x(q) = (x_n)_{n \in \N} \in Num_F(Y)(q)\subset \mcc(X,Y)$ for $q \in Q,$ $$\lim_{n \rightarrow +\infty }F(x_n,q) = 0.$$}
 We call the image space  {{}$$\mathcal{S}_Y(F)= \left\{ \lim_{n \rightarrow +\infty} x \in C^\infty(Q,Y) \, |\, x \in Num_F(Y)\right\}$$} the space of {{} $Q-$parametrized} solutions of (\ref{eq}) with respect to $Num_F(Y).$ \end{Definition}

{{}
	\begin{rem}
		In the definition of the space $\mathcal{S}_Y(F)$, we consider the image of $Num_F(Y)$ with respect to the limit map. This means that, for a fixed parameter $q \in Q$, the space of $Y-$solutions to $F(.,q)=0$  is $\mathcal{S}_Y(F)(q).$ 
	\end{rem}}
If the solution is unique in $Y,$ {{} i.e. if $\mathcal{S}_Y(F)$ has only one element,} so called well-posedness of the solution with respect to the set of parameters $Q$ is ensured by the {{}existence} of smooth numerical schemes. 
The difference that we formalize here between the base space $X$ which serves as a domain for the functions $F(.,q),$ with $q \in Q,$ and the space $Y$ were  will take place what we can call weak solutions in $Y$, is a classical feature in solving functional equations. One may think that this differentiation is due to methods that are not efficient enough to solve the equation. This can be true in some cases. In some other approaches one can get a solution which is a priori in $Y$, but which is proved to be in $X$ by a refined analysis. Let us suggest by the following example (which could be considered as a ``toy example'', {{}but \cite{CPR2016} and references there in show that it is fully interesting}) that the necessity of $Y$ relies directly on the lack of relatively compact neighborhoods in the topology of $X.$ 
   
   \begin{example}
   	Let $X=C^\infty(S^1,\R),$ let $Z = \R$ and let $\Delta = -\frac{d^2}{dx^2}$ be the (positive) Laplacian on $S^1.$ Let us consider the heat operator $e^{-t\Delta}$ which is a smoothing operator for $t>0$ and which converges weakly to $Id_{L^2}$ when $t \rightarrow 0^+.$ Let $w \in L^2(S^1,\R).$
   	We consider the following equation: 
   	$$ F(u,t,w) = 0$$
   	with
   	$u \in X,$
   	$Q = \R_+ \times L^2(S^1,\R)$ and
   	\begin{eqnarray}
   	F: & X \times Q \rightarrow &\R \\
   	& \left(u,(t,w)\right) \mapsto & \frac{1}{2\pi} \int_{S^1} \left((e^{-t\Delta}u) - w\right)^2
   	\end{eqnarray}
   	 $w \in L^2(S^1,\R)$ This equation has an unique solution $u = e^{t\Delta}w$ only when $w \in \{e^{-t\Delta}v \, | \, v \in L^{2}(S^1,\R)\}$ but one can  implement a
   	 numerical scheme even when $w$ has not enough regularity. For example, since $F$ is a convex functional, one may consider a gradient method {{} in order to approximate solutions of the problem}. $$D_uF(v) = \frac{1}{2\pi}\left(\int_{S^1} (e^{-2t\Delta}u)v - 2\int_{S^1}(e^{-t\Delta}w)v\right), $$
   	 {{} hence we get } the $L^2-$gradient 
   	 $$\nabla F(u) = e^{-2t\Delta}u - 2e^{-t\Delta}w$$ {{}Applying the gradient method, we define a sequence $(u_n),$ 
   	 fixing $u_0 \in X,$  setting by induction, 
   	 for $n \in \N$  $u_{n+1} = u_n - \gamma \nabla F(u_n)$ where the constant $\gamma$ is such that $F(u_n) - \gamma||\nabla F(u_n||_{L^2} = 0.$} For ``bad'' choices of $u_0,$ of the parameter $w$ and for $t>0,$ convergence of the sequence $(u_n)$ is a priori accomplished neither for convergence in $X$ nor for  weak $L^2-$convergence, while the sequence remains in the set defined by: $$U =\{ u \in X \, | \, F(u) \leq F(u_0)\},$$ and requires an adapted space $Y$ for which $U$ is bounded and with compact closure. 
   	\end{example}
  
 Since smooth dependence on the parameters is ensured, a notion of (smooth) symmetries can be derived by extension of the classical notions of symmetries. 
 
 One can equip  $\mathcal{S}_Y(F)$ with one of the following diffeologies which appear natural to us:
 
 \begin{itemize}
 \item The push-forward diffeology $\p^{(1)}= \lim_* \p_\C$ which can differ from the subset diffeology inherited from {{}$\mcc(Q,Y).$}
 \item $\p^{(2)},$ the Fr\"olicher completion of $ \p^{(1)},$ {{}i.e. the diffeology $\p^{(2)}$ is the (nebulae) reflexive diffeology associated to the Fr\"olicher structure $$\left(\mathcal{S}_Y(F), \F_{\p^{(1)}}, \mcc_{\p^{(1)}}\right) $$ defined along the lines of Theorem \ref{compl-fro}.}
 \item $\p^{(3)}= \p^{(2)}(Diff),$ where the {{} diffeology $\p^{(2)}(Diff)$ is defined along the lines of Proposition \ref{prop:pdiff}.}
 \end{itemize}323xxb
 {{}
 \subsection{On the way to symmetries}
 We wish to propose in this section a perspective on which objects can be assimilated to the classical notion of groups of symmetries of functional equations. The key problem remains in considering the $Y-$convergence and its smootness as the central aspect of the notion of numerical scheme. After this remark, a symmetry needs to transform one solution to another. This leads to the following two definitions.}
 \begin{Definition}
  The set of $Y-$symmetries, {{} denote}d by $Sym(F,X,Y)$ is defined as the set of smooth maps
   $$\Phi : \mathcal{S}_Y(F)\rightarrow \mathcal{S}_Y(F)  $$  which have a smooth inverse.\end{Definition}
{{}
\begin{Definition}
	The set of $(X,Y)-$sequential symmetries, {{} denote}d by $SSym(F,X,Y),$ is defined as the set of smooth maps
	$$\Phi : Num_F(Y) \rightarrow Num_F(Y)  $$  which have a smooth inverse.
\end{Definition}
 In other words, $$Sym(F,X,Y) = Diff(\mathcal{S}_Y(F)),$$ and $$SSym(F,X,Y) = Diff(Num_Y(F)).$$ We can apply the settings described in section \ref{diff}.

   \begin{Proposition}
   	$Sym(F,X,Y)$  and $SSym(F,X,Y)$ are diffeological groups.
   \end{Proposition}  
\begin{rem}Let $q \in Q.$ Due to Proposition \ref{functf}, here is a commuting diagram of smooth maps
	$$\begin{array}{ccc} Num_F(X,Y) & \rightarrow & S_Y(F) \\
	\downarrow & & \downarrow \\
	Num_F(X,Y)(q)& \rightarrow & S_Y(F)(q) \end{array}$$
	where the horizontal arrows are limit maps, the vertical arrows are evaluation maps at $q.$
	We have here two fiber pseudo-bundles, namely $(Num_F(X,Y), S_Y(F),\lim)$ and $(Num_F(X,Y)(q), S_Y(F)(q),\lim),$ where fibers are spaces of Cauchy sequences. These fibers are not a priori isomorphic. Applying now this pseudo-bundle structure to symmetry groups, we remark successively that 
	\begin{itemize}
		\item $SSym(F,X,Y)$ does not leave invariant the fibers of $(Num_F(X,Y), S_Y(F),\lim_{n \rightarrow +\infty })$
		\item We can define a subgroup of $SSym(F,X,Y)$ by $$Aut(Num_F(X,Y), S_Y(F),\lim_{n \rightarrow +\infty }) $$ $$ = \left\{ g \in Num_F(X,Y) \, | \,  \exists h \in S_Y(F), \, \lim_{n \rightarrow +\infty } \circ g = h \circ \lim_{n \rightarrow +\infty } \right\} $$
		which is analogous to the automorphism group of a (classical, locally trivial and with typical fiber) fiber bundle.
	\end{itemize}    
\end{rem}}
   
   One can define the same way infinitesimal symmetries and sequential infinitesimal symmetries {{} by considering the tangent space at identity of the groups $Sym(F,X,Y)$ and $SSym(F,X,Y).$ These spaces, which are diffeological vector spaces but not a priori diffeological Lie algebras,  can be viewed as sets of smooth sections of the diff-tangent spaces  $^{d}T\mathcal{S}_Y(F)$ and $^{d}TWSym(F,X,Y)$ (see Definition \ref{def:dtangent}).
   
   \begin{Definition}
   	An infinitesimal $Y-$symmetry is defined as an element $u \in {}^iTSym(F,X,Y),$ and an infinitesimal sequential $Y-$symmetry is an element $v \in {}^iTSSym(F,X,Y).$   
   \end{Definition}
There is here an obvious, and yet formal, correspondence with (classical) infinitesimal symmetries viewed as vector fields of the phase space, i.e. formally as elements of the Lie algebra of the group of diffeomorphisms of the phase space \cite{Olv}.Our groups $SSym(F,X,Y)$ and $Sym(F,X,Y)$ appear as some kind of ``maximal'' groups of transformations that enables to deduce an $Y-$solution from another.
   These groups are rather theoretical and cannot be realized easily, even when $X$ is  a space of functions, e.g. when $X$ is of the type $C^\infty(M,V)$ where $M$ is a smooth paracompact manifold and when $V$ is a Euclidian space. In this setting, which is usual in the field of partial differential equations, the approach described in e.g. \cite{Olv} provides a restricted choice of possible symmetries, see e.g. \cite{Vin2013}.  Provided $X =  C^\infty(M,\R),$ $F$ is a smooth map on a finite order jet space over $X$ and $Q = \{q\}$ we can give the following constructions which fix correspondence between our abstract symmetry groups with   the usual ones. The group $G \subset Diff(M)$ of classical symmetries is the group acting on the right on $X$ that leaves invariant $S_X(F).$Infinitesimal (projectable) symmetries are the vector fields that lie in the Lie algebra of $G.$   
   
One can see here that, since one often reaches only the connected cmponent of the identity of $G$, by integrating the Lie algebra of infinitesimal symmetries \cite{Olv} through the analysis of the exponential map (see e.g. \cite{Rob} for one example), the full group $G$ is rarely described.

 Moreover, when $Y \neq X$ and when the set of parameters $Q$ is not reduced to a pôint, $S_Y(F)$ is not known. Let us now describe a well-known example, where the set $S_Y(F)$ should be much bigger than the set $S_X(F),$ and where $Y-$symmetries remain unknown. }
\begin{example} \textbf{Open problem on an example: paradoxical solutions of perfect fluid dynamics}

Let $n\geq 2.$ Let us consider the perfect fluid equations on a compact interval $I=[0;T]$ (at finite time), with variables $u\in C^\infty(I \times \R^n,\R^n)$ (the velocity) and $p \in C^\infty(I \times \R^n, \R),$ and  with (fixed) external force $f$:

\begin{equation} \label{euler}
\left\{ \begin{array}{rc}
\frac{\partial u}{\partial t } - \nabla_v v + \nabla p & = f \\
div (u) & = 0
\end{array}\right.
\end{equation}
Since Scheffer and Shnilerman \cite{Sch,Sh1,Sh2} it is known that these equations have surprising weak solutions, compactly supported in space and time. {{} More precisely, Euler equations have only one (trivial) smooth solution with null initial value, but many weak solutions for a given initial value.}
Advances in this problem have been recently performed in series of papers initiated by \cite{dLS1,dLS2,dLS3}, see also \cite{Vil} for an overview of the two first references. Let us focuse with the following central theorem, see e.g. \cite[Theorem 1.4]{Vil}:

\begin{Theorem}
	Let $\Omega$ be an open subset of $\R^n,$ let $e$ be an uniformly continuous function $ ]0;T[ \times \Omega \rightarrow  \R_+^*$ with $e \in L^\infty\left( ]0;T[; L^1(\Omega) \right).$ $\forall \eta >0,$ there exists a weak solution $(u,p)$ of (\ref{euler}) with $f =0$ such that:
	\begin{enumerate}
		\item $u \in C^0([0;T], L^2_w(\R))^n$ (the subscript ``$w$'' refers to the weak topology)
		\item $u(t,x)=0$ if $(t,x)\in ]0;T[ \times \Omega,$ in particular, $\forall x \in \R^n,$ $u(0,x) = u(T,x) = 0.$
		\item $\frac{|u(t,x)|^2}{2} = \frac{n}{2}p(t,x)= e(t,x)$ almost everywhere on $\Omega$ and forall  $t \in ]0;T[$ 
		\item $\sup_{t \in I}||u(t,.)||_{H^{-1}(\Omega)}<\eta.$
		\item there exists a sequence $(u_n,p_n)$ such that $(u_n, p_n)$ converges to $(u,p)$ for strong $L^2(dt,dx)$-convergence, where each $(u_m,p_m)$ is a (classical) solution of (\ref{euler}) with force $f_m,$ and such that $(f_n)$ converges to $0$ in the sense of distributions.
	\end{enumerate}
	
\end{Theorem}

Actually, the sequence $f_m$ can be obtained from convex integration, the sequence $(u_m,p_m)$ is a Cauchy sequence of (smooth) classical solutions of (\ref{euler}), which is $L^2-$Cauchy, and a natural question is the following:

\vskip 12pt
\noindent
\textbf{Problem} The control function $e$ and the sequence $(f_m)$ appear as parameters to define the Cauchy sequence $(u_m,p_m)_{m \in \N}.$ One can wonder if there is a diffeology on the parameter set $Q= \{(e,(f_m))\} \subset L^\infty\left( ]0;T[; L^1(\Omega) \right)  \times {{}\mcc \left(\mathcal{D}(]0;T[ \times \Omega)\right)}$ (different from the discrete diffeology!) such that the map $(e,(f_m)) \mapsto (u_m,p_m) \in \mcc(\mathcal{D}(]0;T[\times \Omega)^{n+1})$ is smooth. under these conditions, the description of (projectable) symmetries of the set of weak solutions {{}$S_Y(F)$, with $Y = L^\infty([0;T],H^{-1}(\Omega))$ and $F = \frac{\partial u}{\partial t } - \nabla_v v + \nabla p$,} becomes an important open question, while the infinitesimal symmetries \cite{Olv} and the projectable symmetries \cite{Rob} of the set of (non-weak) solutions {{} ( i.e. on the space $S_X(F),$ with $X = C^\infty(\Omega,\R^n)$)} are well-known.

\vskip 12pt This opens the question of the diffeologies envolved in convex integration and in the h-principle.
 \end{example}  
     \section{Implicit functions theorem from the viewpoint of numerical schemes} \label{4}

    We set the following notations, {{}from the standard reference \cite{Om} and along the lines of the recent work \cite{Ma2020-1}}:
    Let  $\hbox{\bf E} = (E_i)_{i \in \N}$ and  
    $\hbox{\bf F} = (F_i)_{i \in \N}$ be two ILB vector spaces, i.e. 
    $\forall i > j, E_i \subset E_j$ and $F_i \subset F_j,$ with {{} smooth inclusion.} {{} We do not assume here, $\forall i > j,$ the image $E_i$ to be dense in $E_j$ for the topology of $E_j,$ and we fit with the assumptions of \cite{Ma2020-1}.} 
    Let $O_0$ be an open neighborhood of $(0{{},}0)$ in $E_0 \times F_0,$ let
    $\hbox{\bf O} = (O_i)_{i \in I}$ with $O_i = O_0 \cap ({E_i \times F_i})$, for $i \in \N \cup \{\infty\}.$ 
    
    {{} Let us now propose a diffeological approach to the main result of \cite{Ma2020-1} that we recall here. For this, we consider a function $f_0$ of class $C^\infty$such that
    	\begin{enumerate}
    		\item $f_0(0; 0) = 0$
    		\item $D_2f_0(0; 0) = Id_{F_0}$.
    	\end{enumerate} Moreover, let us assume that $f_0$ restricts to $C^\infty-$maps
    	$ f_i : U_i \times V_i \rightarrow F_i.$
    	Let $$E_\infty = \underleftarrow{\lim}
    	\left\{ E_i; i \in \N \right\},$$ let $$F_\infty = \underleftarrow{\lim}
    	\left\{ F_i; i \in \N \right\},$$ let $$U_\infty = E_\infty \cap U_0 \quad \hbox{ and } \quad V_\infty = V_0 \cap F_\infty.$$  Finally, let $f_\infty = \underleftarrow{\lim} f_i.$ While we do not assume any other assumption, contrasting with e.g. the classical Nash-Mser theorem \cite{Ham} where additional norm estimates are necessary. Under our weakened conditions, one can state:
 \begin{Theorem} \cite[Theorem 2.2]{Ma2020-1}\label{1.6} There exists a non-empty domain $D_\infty \subset U_\infty,$ possibly non-open
 	in $U_\infty,$ and a function $$u_\infty : D_\infty \rightarrow V_\infty$$ such that
 	$$\forall x \in D_\infty, \quad f_\infty(x; u_\infty(x)) = 0.$$
 	{{} Moreover, there exists a sequence $(c_i)_{i \in \N} \in (\R_+^*)^\N$ and a Banach space $B_{f_\infty}$ such that
 		\begin{itemize}
 			\item $B_{f_\infty} \subset E_\infty$ (as a subset)
 			\item the canonical inclusion map $B_{f_\infty} \hookrightarrow E_\infty$ is continuous
 		\end{itemize} 
 		which is the domain of the following norm (and endowed with it): $$||x||_{f_\infty} = \sup
 		\left\{ \frac{||x_i||}{c_i}| i \in\N \right\}. $$ Then $D_\infty$ contains $\mathcal{B},$ the unit ball (of radius $1$ centered at 0) of  $B_{f_\infty}.$} \end{Theorem}   
In \cite{Ma2020-1}, the question of the regularity of the implicit function is left open, because the domain $D_\infty$ is not a priori open in $O_\infty.$ Moreover, the prosence of the Banach space $B_{f_\infty}$ suggests that the properties of the implicit function $u_\infty$ may depend on the properties of   the function $f_\infty$ under consideration. This lack of regularity induces a critical breakdown in generalizing the classical proof of the Frobenius theorem to this setting. We fill this gap in the sequel, by completing the proof of Theorem \ref{1.6} from \cite{Ma2020-1}  using the Cauchy diffeology,   }
    under the light of numerical schemes. Each of these aspects will be discussed in detailed comments and remarks after the statement and the proof of the following implicit functions theorem. 
    \begin{Theorem} \label{IFTh}
    	{{}Let $$f_i: O_i \rightarrow F_i, \quad i \in \N \cup \{\infty \}$$ be a family of maps, 
    	let $u_\infty$ the implicit function defined on the domain $D_\infty$,  as in Theroem \ref{1.6}.
    	Then, there exists a domain $D$ such that $\mathcal{B} \subset D \subset D_\infty$ such that the function $u_\infty$  is smooth for the subset diffeology of $D.$}
    \end{Theorem}

    \begin{proof} {{} We follow and complete the proof of \cite{Ma2020-1}.} 
    Let $g_i = Id_{F_i} - f_i, $ for $i \in \N \cup \{\infty\}.$
    With this condition, $$D_2g_i(0{{},}0) = 0.$$
    We {{} denote} by $\phi_{i,x}= g_i(x,.).$
    Let $D'_i \subset U_i  \subset E_i$ be the set defined by the following:
    $$ x \in D'_i \Leftrightarrow \left\{ \begin{array}{l} \forall n \in \N,  \phi^n_{i,x}(0)\in O_i \\
    (\phi^n_{i,x}(0))_{n \in \N} \in \mcc(O_i,O_i) \end{array}\right.$$
    {{} where $$\phi^n_{i,x}(0) = \underbrace{\phi_{i,x}\circ ... \circ \phi_{i,x}}_{ n \hbox{ times}}(0).$$}
    We define, for $x \in D'_i,$ $$u_i(x) = \lim_{n \rightarrow +\infty} \phi^n_{i,x}(0) \in \bar{O_i}. $$ For each $x \in D'_i,$ if $u_i(x) \in O_i,$ 
    we get $g_i(x,u_i(x)) = x$ and hence $f_i\left(x,u_i(x) \right) =0.$ {{}By the classical implicit functions theorem on Banach spaces, the domain $D_i$ contains the open ball of $B_i$ centered at $0$ and with radius $c_i/2.$ Moreover, since $f_i$ is smooth, $u_i$ is smooth on $D_i$ and the map 
    $$ x \in D_i \mapsto (\phi^n_{i,x}(0))_{n \in \N}$$ is smooth (when we equip the set of Cauchy sequences in $O_i$ with its Cauchy diffeology).} 
    The family of maps $\{g_i; i \in \N \}$ restricts to a map $g$ on $O_\infty,$ and hence each map $\phi_{i,x}$ restricts the same way to a map $\phi_x.$
    We set {{}$$D= \bigcap_{i \in \N}  D'_i $$
    and  the map $$x \mapsto (\phi^n_{x}(0))_{n \in \N}\in\bigcap_{i \in \N}\mcc(O_i,O_i)=\mcc(O_\infty,O_\infty)$$ is trivially smooth for the Cauchy diffeology on $\mcc(O_\infty,O_\infty).$
     Then, the map $u_\infty$ which is the restriction of and map $u_i$ to $D$ is a smooth map. Finally, $D$ contains element $x \in O_\infty$ such that $$\forall i \in \N, \quad ||x||_i < c_i/2$$
 which is equivalent to the condition $$\sup_{i \in \N } \frac{||x||_i}{c_i}\leq 1/2 < 1,$$ which ends the proof.} \end{proof}
    After this proof, adapted from the classical proof of the implicit functions theorem in Banach spaces, see e.g. \cite{Pen},  we must notice several points:
    
    \begin{enumerate}
    	
    	\item Passing to the projective limit, the question of the nature of the domain $D$ can be adressed. As an infinite intersection of open sets, it cannot be stated that $D$ is an open subset of $E_\infty.$ We could add the remark that, skeptically, $D$ can be restricted to $\{(0{{},}0)\}$ {{} but the presence of the open ball of radius 1 of $B_{f_\infty}$ shows that there exists  ways to modify the topology of $O_\infty$ in order to make the domain  $D$ open, and one of these ways is described explicitely.} Thinking with an intuition based on finite dimensional problems, one could think that a domain $D$ which is not a neighborhood of $0$ is not interesting, because all norms are equivalent. This is the motivation for Nash-Moser estimates \cite{Ham}, with which one manages to show that the domain $D$ is an open neighborhood of $0.$ This comes from the natural desire to control the "size"
    	of the domain $D.$ But noticing that $E_1 \subset E_0$ is not an open subset of $E_0$ unless $E_1 = E_0,$ one can see that $D$ is not open {{} for the topology of $O_\infty$} is not disqualifying. 
    	
    	\item Since Nash-Moser estimates appears rather artificial even if very useful, and only designed to "pass to the projective limit" some open domain, and noticing that the condition of contraction enables to bound the difference $||\phi^{n+1}_{i,x}(0) -  \phi^{n}_{i,x}(0)||_{F_i}$ by a converging geometric sequence, one can try to generalize this procedure and apply the notion of bornology. This gives the theorem stated in \cite{HN1971}. But Theorem \ref{1.6} deals with a {{} different} class of functions. 
    	
    	\item Summarizing these two facts and the detailed proof given before, one gets smoothness of the solution $u$ of the equation $f(x,u(x))=0$ and its uniqueness on $D \cap O_0.$ This relates implicit functions theorem to well-posedness, i.e. the problem of unique existence of a unique solution smoothly dependent on the set of parameters $Q = D.$ The same proof, applied to the functions $f_i$ on the domain $O_\infty,$ produce weak solutions $u_i(x) \in F_i,$ depending smoothly on the parameter $x \in D,$ to the equation $f_\infty(x,u(x))=0.$ 
    	
    	\item Pushing forward this last point, again applying the classical implicit functions theorem on Banach spaces \cite{Pen}, one gets implicit functions $u_i : D_i \rightarrow F_i$ where $D_i$ is a neighborhood of $0$ in $O_i.$ Restricting to $Q_i = D_i \cap U_\infty,$ one can say that $Q_i$ is a set of parameters for solvig the equation $F(u,q)=0$ with $ F(u,x)= f_\infty(x,u(x)).$ Each function $u_i$ appears as the unique element of $\mathcal{S}_{F_i}(F),$ with admissible set of parameters $Q_i.$The set $D$ appears as the set of admissible parameters at each index $i \in \N ,$ i.e. $D  = \cap_{i \in \N} Q_i.$
    \end{enumerate}   
    
    \vskip 12pt
    The same way, we can state the corresponding Frobenius theorem:
    
    \begin{Theorem}\label{lFrob}
    	
    	Let 
    	$$ f_i : O_i \rightarrow L(E_i,F_i), \quad i \in {{}\N}$$ 
    	be a collection of smooth maps satisfying the following condition: 
    	$$ i > j \Rightarrow f_j|_{O_i} = f_i$$ and such that, $$\forall (x,y) \in O_i, \forall a,b \in E_i$$
    	$$(D_1f_i(x,y)(a)(b) + (D_2f_i(x,y))(f_i(x,y)(a))(b) =$$
    	$$(D_1f_i(x,y)(b)(a) + (D_2f_i(x,y))(f_i(x,y)(b))(a) .$$

    	Then,
    	$\forall (x_0, y_0) \in O_{\infty}$, there exists a diffeological subspace  $ D $ of $O_\infty$ that contains $(x_0, y_0)$ and a smooth map
    	$J : D \rightarrow  F_\infty$
    	such that
    	$$ \forall (x,y) \in D, \quad D_1J(x,y) = f_i(x, J(x,y)) $$
    	and, if {{}$D_{x_0}$ is the connected component of $(x_0,y_0)$ in $\{(x,y) \in D \, | \, x = x_0  \},$ }
    	$$J_i(x_0,.) = Id_{D_{x_0}}.
    	$$ 
    	{{} Moreover, there exists a sequence $(c_i)_{i \in \N} \in (\R_+^*)^\N$ and a Banach space $B_{f_\infty}$ such that
    		\begin{itemize}
    			\item $B_{f_\infty} \subset E_\infty\times F_\infty$ (as a subset)
    			\item the canonical inclusion map $B_{f_\infty} \hookrightarrow E_\infty\times F_\infty$ is continuous
    		\end{itemize} 
    		which is the domain of the following norm (and endowed with it): $$||x||_{f_\infty} = \sup
    		\left\{ \frac{||x||_{E_i \times F_i}}{c_i}| i \in\N \right\}. $$ Then $D_\infty$ contains $\mathcal{B},$ the unit ball (of radius $1$ centered at 0) of  $B_{f_\infty}.$}
    	
    \end{Theorem}

    \begin{proof}
    We consider  
    $$G_i = C^1_0([0,1],F_i) = \{ \gamma \in C^1([0,1],F_i) | \gamma(0)=0 \}$$ and 
    $$ H_i = C^0([0,1],F_i),$$
    endowed with their usual topologies, and nebulae underlying diffeology. Obviously, if $ i < j$, 
    the injections $ G_j \subset G_i$ and $ H_j \subset H_i$ are smooth.
    
    Let us consider open subsets $B_0\subset E_0$ and   $B'_0\subset F_0,$ 
    such that $$(x_0,y_0) \subset B_0 \times B'_0 \subset O_0.$$ Let us consider the open set   
    $$B''_0 = \{ \gamma \in G_0 | \forall t \in [0;1], \gamma(t) \in B'_0\}.$$  We set 
    $ B_i = B_0 \cap E_i$,  $ B_i' = B_0' \cap F_i$ and 
    $ B_i'' = B_0'' \cap G_i$. Then, we define, for $ i \in \N\cup \{\infty\}$,  
    $$ g_i : B_i \times B'_i \times B''_i \rightarrow H_i $$
    $$ g(x,y,\gamma)(t) = \frac{d \gamma }{ dt}(t) - f_i(t(x-x_0) + x_0, y + \gamma(t)).(x - x_0).$$
    In order to apply the last implicit function theorem, we must calculate, avoiding the subscripts $i$ for easier reading: 
    $$ D_3g(x_0,y_0,0)(\delta) = \frac{d \delta }{ dt} .$$
    Thus, we can apply Theorem \ref{IFTh} to $$f = \left(\int_0^{(.)}\right)\circ g .$$

    We can define the function $\alpha$ as the unique function on the domain $D \subset B\infty \times B'_\infty$ such that
    
    $$ \left\{ \begin{array}{l}
    \alpha(x_0,y_0) = 0 \\
    g_\infty(x,y,\alpha(x,y)) = 0 , \quad \forall (x,y) \in D \\
    \end{array}
    \right. .$$
    We set $J(x,y) = y + \alpha(x,y)(1).$ {{} The constants $c_i$ are obtained through the application of Theorem \ref{IFTh} }  \end{proof}
 \section{The Fr\"olicher space of triangulations and the finite elements method} \label{5}
 Let us consider the 2-dimensional Dirichlet problem. Let $\Omega$ be a bounded connected open subset of $\R^2,$ and assume that the border $\partial \Omega = \bar{\Omega} - \Omega$ is a polygonal curve.  The Dirichlet problem is given by the following PDE:
 
 $$ \left\{ \begin{array}{c}
 \Delta u = f \\
 u|_{\partial \Omega} = 0
 \end{array}
   \right. $$
   where {{}$f$ is a smooth, compactly supported function in $\Omega$ (i.e. $f  \in \D(\Omega, \R)),$} $\Delta$ is the Laplacian, and $u$ is the solution of the Dirichlet problem. 
   
   \vskip 12pt
   Let us analyze the set of triangulation as a smooth set of parameters for the finite element method, with as an example the Dirichlet problem. 
   \subsection{{{}A quick summary} of the finite elements method (degree 1) for the Dirichlet problem.}
   
   One classical way to solve the problem is to approximate $u$ by a sequence $(u_n)_{n \in \N}$ in the Sobolev space $H^1_0(\Omega,\R)$ which converges to $u$ for the $H^1_0-$convergence. For this, based on a triangulation $\tau_0$ with $0-$vertices $(s^{0}_k)_{k \in K_0},$ where $K_0$ is an adequate set of indexes, and we consider the $H^1_0-$orthogonal family $\left(\delta_{s_k^{(0)}}\right)_{k \in K_0}$ of continuous, piecewise affine maps on each interior domain of triangulation, defined by $$\delta_{s_k^{(0)}} (s_j^{(0)}) = \delta_{j,k} \hbox{ (Kronecker symbol).}$$   
   With this setting, $u_0$ is a linear combination of $\left(\delta_{s_k^{(0)}}\right)_{k \in K_0}$ such that
   $$\forall k \in K_0, \left(\Delta u_0, \delta_{s_k^{(0)}}\right)_{H^{-1}\times H^1_0} = \left( f, \delta_{s_k^{(0)}}\right)_{L_2} .$$
   This is a finite dimensional linear equation, which can be solved by inversion of a $|K_0|-$dimensional matrix. 
   Then we refine the triangulation $\tau_0$ adding the centers of $1-$vertices to make the triangulation $\tau_1,$ and by induction, we get a sequence of triangulation $(\tau_n)_{n \in \N}, $ and by the way a sequence of  families $$ \left( \left(\delta_{s_k^{(n)}}\right)_{k \in K_n} \right)_{n \in \N},$$ which determines the sequence $(u_n)_{n \in \N}$ which converges to $u.$
   Then through operator analysis, we know that $u$ is smooth and that there exists a smooth inverse to the Lapacian $\Delta^{-1}$ adapted to the Dirichlet problem such that $u = \Delta^{-1}f.$ For the rest of the section, we equip $H^1_0(\Omega,\R)$ by the Fr\"olicher structure
   generated by $$\left\{ \left(.,f\right)_{H^1_0} \, | \, f \in C^\infty_c(\Omega,\R)\right\}.$$
   
   \subsection{Differential geometry of the space of triangulations}
   
   Let us now fully develop an approach based on the remrks given in  \cite{Ma2016-2}. For this, the space of triangulations of $\Omega$ is considered itself as a Fr\"olicher space, and the mesh of triangulations which makes the finite element method converge will take place, as the function $f,$  among the set of parameters $Q.$ We describe here step by step the Fr\"olicher structure on the space of triangulations.  
   By the way, we begin with a lemma which is adapted from so-called gluing results present in \cite{Nt2002,pervova2017,pervova2018} to the context which is of interest for us.
   \begin{Lemma} \label{cov}
   	Let us assume that $X$ is a topological space, and that there is a collection $\{(X_i,\F_i,\mcc_i)\}_{i \in I}$ of Fr\"olicher spaces, 
   	together with continuous maps $\phi_i: X_i \rightarrow X.$
   	Then we can define a Fr\"olicher structures on $X$
   	setting $$\F_{I,0} = \{f \in C^0(X,\R) | \forall i \in I, \quad f \circ \phi_i \circ \mcc_i \subset C^\infty(\R,\R)\},$$ wa define $\mcc_I$ the contours generated by the family $\F_{I,0},$ and  $\F_I = \F(\mcc_I).$
   \end{Lemma}
   
   Let $M$ be a smooth manifold for dimension $n.$ Let \begin{equation} \label{embtriangulation}\Delta_n= \{(x_0,...,x_n)\in \R_+^{n+1} | x_0 +... +x_n = 1\}\end{equation} be the standard  $n-$ simplex, equipped with its subset diffeology. It is easy to show that this diffeology is reflexive through Boman's theorem already mentionned, and hence we can call it Fr\"olicher space $(\Delta_n, \F_{\Delta_n}, \C_{\Delta_n}),$ and we {{} denote} its associated reflexive diffeology by $\p(\Delta_n).$ 
   
   \begin{Definition} \label{smtriang} A \textbf{smooth triangulation} of $M$ is a family $\tau = (\tau_i)_{i \in I}$
   where $I \subset \N$ is a set of indexes, finite or infinite, each $\tau_i$ is a smooth map $\Delta_n \rightarrow M,$ and such that:
   \begin{enumerate}
   	\item $\forall i \in I, \tau_i$ is a (smooth) embedding, i.e. a smooth injective map such that {{}$(\tau_i)_*\left(\p({\Delta_n})\right)$} is also the subset diffeology of $\tau_i(\Delta_n)$ as a subset of $M.$
   	\item $\bigcup_{i \in I }\tau_i(\Delta_n) = M.$ (covering)
   	\item $\forall (i,j) \in I^2,$ $\tau_i(\Delta_n) \cap \tau_j(\Delta_n) \subset  \tau_i(\partial \Delta_n) \cap \tau_j(\partial \Delta_n).$ (intersection along the borders)
   	\item  $\forall (i,j) \in I^2$ such that $ D_{i,j}=\tau_i(\Delta_n) \cap \tau_j(\Delta_n) \neq \emptyset,$ for each $(n-1)$-face $F$ of $D_{i,j},$ the ``transition maps" $``\tau_j^{-1} \circ \tau_i'' : \tau_i^{-1}(F) \rightarrow \tau_j^{-1}(F)$ are affine maps.  
   \end{enumerate} 
\end{Definition} 
   Under these  conditions, we equip the triangulated manifold $(M,\tau)$ with a Fr\"olicher structure $(\F_I,\mcc_I),$ generated by the smooth maps $\tau_i$ applying Lemma \ref{cov}. The following result is obtained from the construction of $\F$ and $\mcc:$
   \begin{Theorem}
   	The inclusion $ (M,\F,\mcc) \rightarrow M$ is smooth.
   \end{Theorem}
   
   \begin{proof} Here the manifold $M$ is considered as a reflexive diffeological space equipped with its reflexive diffeology $\p_\infty(M)$ and with its associated Fr\"olicher structure $(\F(\p_\infty(M)), \C(\p_\infty(M))).$ 
   	Let $f \in C^\infty((M,\p_\infty(M)),\R),$ i.e. a (classical) smooth map 
   	$f \in C^\infty(M,\R).$ Since each $\tau_i$ is a smooth map $\Delta_n \rightarrow M,$ 
   	$$f \circ \tau_i \circ \C \in C^\infty(\R,\R)$$ and hence $$C^\infty((M,\p_\infty(M)),\R) \subset \F_{I,0} \subset \F_I.$$
   	
   	\end{proof}
   
   \begin{rem} Maps in $\F_I$ can be intuitively identified as some piecewise smooth maps $M \rightarrow \R,$ which are of class $C^0$ along the 1-skeleton of the triangulation.
   	We have proved also that $\mcc_I \subset \p_\infty(M).$ Some characteristic elements of $\mcc_I$ can be understood as paths which are smooth (in the classical sense) on the interiors of the domains of the simplexes of the triangulation, and that fulfill some more restrictive conditions while crossing the 1-skeleton of the triangulation. For example, paths that are (locally) stationnary at the 1-skeleton are in $\mcc_I.$
   \end{rem}

\begin{rem}
	While trying to define a Fr\"olicher structure from a triangulation, one could also consider $$\C_{I,0} = \left\{ \gamma \in C^{0}(\R,M)\, | \, \forall i \in I, \forall f \in C^\infty_c(\phi_i(\Delta_n),\R), f \circ \gamma \in C^\infty(\R,\R) \right\}$$ where $C^\infty_c(\phi_i(\Delta_n),\R)$ stands for compactly supported smooth functions $M \rightarrow \R$ with support in $\phi_i(\Delta_n).$ Then define $$\F_I' = \left\{f : M \rightarrow \R \, | \, f \circ \C_{I,0} \in C^\infty(\R,\R)\right\}$$
	and $$\C_I' = \left\{C : \R \rightarrow M \, | \, \F_I' \circ c \in C^\infty(\R,\R)\right\}.$$
	 We get here another construction, but which does not understand as smooth maps $M \rightarrow \R$ the maps $\delta_k$ already mentionned.
\end{rem}
   Now, let us fix the set of indexes $I$ and fix a so-called \textbf{model triangulation} $\tau.$ {{} This terminology is justified by two ideas: 
   	\begin{itemize} 
   		\item Anticipating next constructions, this model triangulation $\tau$ will serve at defining a sequence of refined trinagulations. This is our ``starting triangulation'' for the refinement procedure in the finite elements method. 
   		\item Changing $\tau$ into $g \circ \tau,$ where $g$ is a diffeomorphism, we get another model triangulation, which has merely the same properties as $\tau.$ But each ``starting'' trinagulation cannot be obtained by transforming a fixed triangulation by using a diffeomorphism. For example, on the 2-sphere, a tetrahedral triangulation $\tau_1$ and an octahedral triangulation $\tau_2$ separately generate two sequences of refined triangulations, and there is a topological obstruction for changing $\tau_1$ into $\tau_2$ by the action of a diffeomorphism of the sphere.     
   		\end{itemize} } We {{} denote} by $\mathcal{T}_\tau$ the set of triangulations $\tau'$ of $M$ such that the corresponding 1-skeletons are diffeomorphic to the 1-skeleton of $\tau$ (in the Fr\"olicher category). {{} The set $\mathcal{T}_\tau$ contains, but is not reduced to, the orbit of $\tau$ by the action of the group of diffeomorphisms. Indeed, one can reparametrize each simplex with adequate compatibility on the border. Intuitively speaking, reparametrizations need not to be smooth in the usual sense while ``crossing the border of a simplex''. This choice is motivated by the Fr\"olicher structure that we identify as useful for the finie elements method, ddefined hereafter.} 
   \begin{Definition} \label{d3}
   	Since $\mathcal{T}_\tau \subset C^\infty(\Delta_n, M)^I,$ we can equip  $\mathcal{T}_\tau$ with the subset Fr\"olicher structure, in other words, the Fr\"olicher structure on $\mathcal{T}_\tau$ whose generating family of contours $\mcc$ are the contours in $C^\infty(\Delta_n, M)^I$ which lie in $\mathcal{T}_\tau.$
   \end{Definition} 
We define the full space of triangulations $\mathcal{T}$ as the disjoint union of the spaces of the type  $\mathcal{T}_\tau,$ with disjoint union Fr\"olicher structure. With this notation, in {{}the} sequel and when it carries no ambiguity, the triangulations in $\mathcal{T}_\tau$ is equipped with a fixed set of indexes $I$ (which is impossible to fix for $\mathcal{T}$). We need now to describe the procedure which intends to refine the triangulation and define a sequence of triangulations $(\tau_n)_{n \in \N}.$ 
We can now consider the refinement operator, which is the operator which divides a simplex $\Delta_n$ into a triangulation. 

\begin{Definition}
	Let $m \in \N,$ with $m \geq 3.$ Let 
	$$\mu = \left\{\mu_i: \Delta_n \rightarrow \Delta_n \, | \, i \in \N_m\ \right\}$$ be a smooth triangulation of $\Delta_n$ 
	Let $\tau \in \mathcal{T}.$ Then we define {{}$$\mu(\tau) = \{f_i \circ \mu_i \, | \,  i \in \N_m \hbox{ and } \tau = (f_i)_{i \in I}\}.$$ We say that $\mu$ defines a \textbf{refinement map} if $\forall n \in \N^*, \mu^n(\tau)$ is a triangulation. }
\end{Definition}

With this definition, $\mu(\tau)$ is trivially a triangulation of $M$ if $\tau$ is a triangulation of $M.$
The conditions imposed in the definition ensures that the refinement map maps {{}a} triangulation to another {{}triangulation}, that is, if $\tau$ is a triangulation, $\mu(\tau)$ is {{}also} a {{}triangulation}. The delicate needed condition is that the new 0-vertices added to $tau$ in $\mu(\tau)$ are matching. 
\begin{Theorem}
	The map $\mu: \mathcal{T} \rightarrow \mathcal{T}$ is smooth.
\end{Theorem}

\begin{proof}
	Composition map $$ C^\infty(\Delta_n, \Delta_n) \times \C^\infty(\Delta_n, M) \rightarrow C^\infty(\Delta_n, M)$$ is smooth, so that it extends canonically (coefficientwise) to a smooth map $$\Phi : C^\infty(\Delta_n, \Delta_n)^{\N_n} \times \C^\infty(\Delta_n, M)^I \rightarrow C^\infty(\Delta_n, M)^{\N_m \times I}.$$
	Let us fix $\mu \in C^\infty(\Delta_n, \Delta_n)^{\N_n}$ a smooth triangulation of $\Delta_n,$ for a fixed model triangulation $\tau = \{\tau_i\}_{i \in I}$ the map $\mu$ is a restriction $\mathcal{T}_\tau \rightarrow \mathcal{T}_{\mu(\tau)}$ of the map $\Phi(\mu,.).$ So that, $\mu \in C^\infty(\mathcal{T}_\tau , \mathcal{T}_{\mu(\tau)})$ and extending this result to $\mathcal{T}$ as a disjoint union,  we get $\mu \in C^\infty(\mathcal{T},\mathcal{T}).$   
	\end{proof}

\begin{Definition}
	Let $\tau \in \mathcal{T}.$ We define the $\mu-$refined sequence of triangulations $\mu^\N(\tau) = (\tau_n)_{n \in \N}$ by $$ \left\{ \begin{array}{ccl} \tau_0 & = & \tau \\ \tau_{n+1} & = & \mu(\tau_n) \end{array} \right.$$ 
\end{Definition}

\begin{Proposition} \label{seqref}
	The map $$\mu^\N : \mathcal{T} \rightarrow \mathcal{T}^\N$$ is smooth (with $\mathcal{T}^\N$ equipped with the infinite product Fr\"olicher structure).
	\end{Proposition}
In the case of the Dirichlet problem, we consider a subspace of $\mathcal{T}_\tau.$
     
     \begin{proof}
     	It follows from smoothness of $\mu:\mathcal{T} \rightarrow \mathcal{T}.$
     	\end{proof}
     
     \begin{Lemma}
     	Let us fix an indexation of the 0-vertices of $\Delta_n.$ Let $\tau = (\tau_i)_{i \in I} \in {{}\mathcal{T}}$ and let$(i,j) \in I \times \N_{n+1}.$ Let $x_j(\tau_i) \in M$  be the image by $\tau_i$ of the $j-$th 0-vertex of $\Delta_n.$ Then for fixed indexes i and j, the map $\tau \in \mathcal{T} \mapsto x_j(\tau_i) \in M$ is smooth.
     \end{Lemma}
 
 \begin{proof} It follows from the smoothness of evaluation maps.
 	\end{proof}
   Let $\Omega$ be a bounded connected open subset of $\R^n,$ and assume that the border $\partial \Omega = \bar{\Omega} - \Omega$ is a polygonal curve. Since $\R^n$ is a vector space, we can consider the space of affine triangulations: 
   $$Aff\mathcal{T}_\tau = \left\{ \tau' \in \mathcal{T}_\tau | \forall i , \tau_i' \hbox{ is (the restriction to } \Delta_n \hbox{ of) an affine map } \right\}.$$
   We define $Aff\mathcal{T}$ from $Aff\mathcal{T}_\tau$ the same way we defined $\mathcal{T}$ from $\mathcal{T}_\tau,$ via disjoint union.
   We equip $Aff(\mathcal{T}_\tau)$ {{}and} $Aff(\mathcal{T})$ with their subset diffeology. We use here the notations of last Lemma.
   \begin{Theorem}
   	Let $$c : \R \rightarrow Aff(\mathcal{T}_\tau)$$ be a path on $Aff(\mathcal{T}_\tau).$ Then $$ c \hbox{ is smooth } \Leftrightarrow \forall (i,j) \in I \times \N_{n+1}, t \mapsto x_j(c(t)_i) \hbox{ is smooth. }$$
   	\end{Theorem}
   	\begin{proof}
   	Let
$x \in \Delta_n.$
We consider the normalized barycentric coordinates $$(\alpha_1(x),...,\alpha_{n+1}(x))$$ 
which correspond to the coordinates of
$x \in \R^{n+1}.$
The map $$ x \in \Delta_n \mapsto (\alpha_1(x),...,\alpha_{n+1}(x))$$ is smooth. Let $c: \R \rightarrow \mathcal{T}$ such that, $\forall    (i,j) \in I \times \N_n, $ the maps $t \mapsto x_j(c(t)_i)$ are smooth. 
We fix $i \in I$ and consider a smooth plot $p \in \p(\Delta_n).$ Then the map $(\alpha_1 \circ p,...,\alpha_{n+1}\circ p)$ is smooth and since $\tau_i$ is affine, {{}$$\tau_i \circ p = \sum_{j = 0}^{n+1} (\alpha_j \circ p). x_j(\tau_i).$$ } We replace $\tau_i$ by $c(t)_i$ in this formula, made of smooth operations, which shows that the maps $ t \mapsto c(t)_i$ are smooth for the diffeology defined in Definition \ref{d3} applying Proposition \ref{functf}.
   	\end{proof}

   	\begin{Proposition}
   	Let {{}$\mu$} be a fixed affine triangulation of $\Delta_n.$
   	The map $\mu^\N$ restricts to a smooth map from the set of affine triangulations of $\Omega$ to se set of sequences of affine triangulations of $\Omega.$
   	\end{Proposition}
   	
   	\begin{proof}
   	Follows from Proposition \ref{seqref}.
   	\end{proof}
   \subsection{Back to the Dirichlet problem}
    With a sequence of affine triangulations $(\tau_n)_{n \in \N}$ defined as before on a suitable domain $\Omega$ of $\R^n,$ we wish to establish smoothness of the family of maps $\delta$ defined before with respect to the underlying triangulation. For this, we extend first the family of $H^1_0-$functions $\delta$ to $\mathcal{T}.$
    
    \begin{Definition}
    Let $\tau \in \mathcal{T},$ indexed by the set $I.$ Let $a$ be a $0-$vertex of $\tau$ We {{} denote} by $St(a)$ the domain described as $$ \cup\{ Im(\tau_i) \, | \, i \in I \hbox{ and } a \in Im(\tau_i)\}.$$ Let  us define the following maps: 
    \begin{itemize}
    \item for $(i,j) \in I \times \N_{n+1}, $
    let $\delta^\tau_{i,j}: \Omega \rightarrow \R$ be the map defined by $$\delta_{i,j}^\tau (x) = \left\{ \begin{array}{ccl} 0 & \hbox{ if } & x \notin Im(\tau_i)  \\
    \alpha_{j}\left( \tau_i^{-1}(x) \right) & \hbox{ if } & x \in Im(\tau_i)
    \end{array}\right.$$
    \item Let $\{ x_k \}_{k \in K}$ be the set of $0-$vertices in $\Omega$ of the triangulation $\tau,$ indexed by $K.$ Let $\delta_{x_k}$ be the map defined by
    $$ \delta^\tau_{x_k}(x) = \left\{ \begin{array}{ccl} 0 & \hbox{ if } & x \notin St(x_k)  \\
     \delta_{i,j}^\tau(x) & \hbox{ if } & x \in Im(\tau_i)\cap St(x_k) \hbox{ and } x_k = x_j(\tau_i)
    \end{array}\right.$$
    \end{itemize}
    
    \end{Definition}
    
    We remark that this definition is consistent by condition (4) of Definition \ref{smtriang}, which ensures that ``gluing along the borders" is possible, that is,  $\forall((i,j),(k,l)) \in (I \times \N_{n+1})^2,$ if $x_k = x_j(\tau_i) = x_l(\tau_k), $ for $x \in Im(\tau_i) \cap Im(\tau_k),$ $$\delta^\tau_{i,j}(x) = \delta^\tau_{k,l}(x).$$
    With the previous notations, we have:
    \begin{Lemma}
    $\forall \tau \in \mathcal{T}, $ $\forall k \in K,$ $\delta_{x_k}^\tau \in H^1_0 \cap C^0(\Omega).$ 
    \end{Lemma}
    
    \begin{proof}
    The map $\delta_{x_k}^\tau$ is: 
    \begin{itemize}
    \item smooth on each interior of domain $\dot{Im(\tau_i)}$
    \item $C^0$ in $\Omega$
    \end{itemize}
    So that, it is a continous map, piecewise smooth. 
    \end{proof}
    
    By the way
    we define a map 
    
    $$\delta: \mathcal{T}_\tau \rightarrow \left(H^1_0 \cap C^0(\Omega)\right)^I$$
    which extends, if $I$ is finite and if $(\tau_n)$ is a $\mu-$refined sequence, to a map
    $$\mu^\N\left(\mathcal{T}_\tau\right) \rightarrow \left(\left(H^1_0 \cap C^0(\Omega)\right)^\infty\right)^\N$$ {{} where $\left(H^1_0 \cap C^0(\Omega)\right)^\infty = \cup_{n \in \N^*} \left(H^1_0 \cap C^0(\Omega)\right)^n,$ with product diffeology,}
    or, if $I = \N,$ to a map $$\mu^\N\left(\mathcal{T}_\tau\right) \rightarrow \left(\left(H^1_0 \cap C^0(\Omega)\right)^\N\right)^\N.$$
    
    \begin{Theorem}
    Let $\tau \in \mathcal{T}.$
    The map $$\delta: \mathcal{T}_\tau \rightarrow \left(H^1_0 \cap C^0(\Omega)\right)^I$$ is smooth. 
    \end{Theorem}
    
    \begin{proof}
    Let us fix $k \in K$ and $f \in C^\infty_c(\Omega,\R).$
    Let $h = \Delta f \in C^\infty_c(\Omega,\R).$ Let $p $ be a plot in the nebulae diffeology of Fr\"olicher structure on $\mathcal{T}_\tau.$ Let $\beta : D(p) \rightarrow H^1_0(\Omega,\R)$ be the map defined by, $$\forall x \in D(p), \quad\gamma(x) = \delta_{x_k}^{p(x)}\in H^1_0(\Omega,\R).$$
    Let $i \in I.$ We define $h_{i,x}: \Delta_n \rightarrow \R$ by $$h_{i,x} = h \circ\gamma(x)_i.$$
    Then 
    \begin{eqnarray*}
    	(\delta_{x_k}^{p(x)}, f)_{H^1_0 } & = & (\delta_{x_k}^{p(x)}, h)_{L^2}\\
    	& = & \sum \int_{\Delta_n} \alpha_j(y) h_{i,x}(y) |J(\tau_i(y))| dy
    	\end{eqnarray*}
     where, in this last equation, the sum $\Sigma$ is among the indexed in $i$ which correspond to $St(x_k),$ $y$ is such that $\tau_i(y)=x$ and $J(\tau_i(y))$ is the Jacobian determinant. Let $c: \R \rightarrow D(p)$ be a smooth path. In order to prove the theorem, via Boman theorem already cited,
     it is sufficient to prove that $ t \mapsto 	(\delta_{x_k}^{p\circ c(t)}, f)_{H^1_0 }$ is smooth for each smooth path $c.$  We have that $$t \mapsto h_{i,c(t)}$$ is smooth in $C^\infty(\Delta_n,\R)$ and by the way, $$t \mapsto \int_{\Delta_n} \alpha_j(y) h_{i,{{}c(t)}}(y) |J(\tau_i(y))| dy$$ is smooth.
    \end{proof}
    
    Now, let us fix $\mu$ a triangulation of $\Delta_n,$ consider the $\mu-$refinement sheme in $\mathcal T$ which introduces, for each $\tau \in \mathcal T,$ a sequence $\tau_n,$ and a family of functions $\delta^{\tau_n}.$
    For fixed index $n$, we solve the problem
    
    $$(\Delta u, \delta_{x_k}^{\tau_n})_{H^{-1}\times H^1_0} = (f,\delta_{x_k}^{\tau_n})_{H^{-1}\times H^1_0}$$
    in the vector space spanned by the family  of functions $\delta_{x_k}^{\tau_n},$ where each $x_k$ is a 0-vertex of $\tau_n.$ If $K$ is the set of indexes $k$ {{} and with cardinal $|K|,$ we get a square matrix $A^{\tau_n}$ with complex coefficients} which is invertible, defined by $$A_{k,l}^{\tau_n} = (\Delta \delta_{x_k}^{\tau_n}, \delta_{x_l}^{\tau_n})_{H^{-1}\times H^1_0},$$ from $v \in \mathbb{C}^{|K|}$ defined by  $$ (f,\delta_{x_k}^{\tau_n})_{H^{-1}\times H^1_0}$$
    we define $$u_n = \left(A^{\tau_n}\right)^{-1} v.$$
    \begin{Theorem}
    	The map $(\tau_0',f ) \in \mathcal{T} \times C^\infty(\Omega,\R)\mapsto (u_n)_{n \in \N}$ is a smooth $H^1_0-$ numerical scheme for the Dirichlet problem.
    \end{Theorem}  
 \begin{proof}
 	{{} The scalar products in the definition of $v$ and of {{}the} coefficients of the matrix $A^{\tau_n}$ are smooth, and the inversion in the group of invertible matrices is smooth too. By the way, the map $$\tau \mapsto u_n$$ is smooth for the infinite product diffeology $\p(H^1_0)_\N.$ In order to get a smooth numerical scheme, we need to be sure that the limit is smooth with respect to $(\tau_0',f ).$ The limit does not depend on $\tau_0$ and is already known to be smooth since the map $f \mapsto \lim_{n \rightarrow +\infty}$ is a well-known pseudo-differential operator of order $-2.$ 
 	}
 	
 \end{proof}


\begin{thebibliography}{99}


\bibitem{BH}
 Baez, J.; Hoffnung, A.; Convenient categories of smooth spaces. \emph{Trans. Amer. Math. Soc.}, \textbf{363}, 5789-5825 (2011).

	\bibitem{BIgKWa2014} Batubenge, A.; Iglesias-Zemmour, P.; Karshon, Y.; Watts, J.A.; Diffeological, Fr\"olicher, and differential spaces (Preprint 2014)
	
	\bibitem{BN2005} Batubenge, A.; Ntumba, P.; On the way to Fr\"olicher Lie groups \textit{Quaestionnes Math.} \textbf{28} (2009) 73-93
	
	\bibitem{BT2014} Batubenge, A.; Tshilombo, M.H.; Topologies on product and coproduct Fr\"olicher spaces; \textit{Demonstratio Math.} \textbf{47}, no4 (2014) 1012-1024
	
	{{}
	\bibitem{Bil1999} Billingsley, P.; {\it Convergence of probability measures, 2nd edition.} Wiley Interscience (1999)}
	
	\bibitem{Bou}  Bourbaki, N.; {\it El\'ements de math\'ematiques}; Masson, Paris (1981)
	
	\bibitem{CPR2016} Carlsson, M.; Prado, H.; Reyes, E. G.
	Differential equations with infinitely many derivatives and the Borel transform. 
	\textit{Ann. Henri Poincar\'e} \textbf{17} no. 8, 2049-2074 (2016). 
	\bibitem{CSW2014}  Christensen, J.D.; Sinnamon, G.; Wu, E.;The D-topology for diffeological spaces, \textit{Pacific Journal of Mathematics}  \textbf{272} no 1 (2014), 87-110
	
	\bibitem{CW2014}  Christensen, J.D.; Wu, E.; Tangent spaces and tangent bundles for diffeological spaces;  {\it Cahiers de Topologie et G\'eom\'etrie Diff\'erentielle},  \textbf{LVII}, 3-50 (2016)
	
	\bibitem{CN} Cherenack, P.; Ntumba, P.; Spaces with differentiable
	structure an application to cosmology \textit{Demonstratio Math.}
	\textbf{34} no 1 (2001), 161-180
	
	\bibitem{Don} Donato, P.; \textit{Rev\^etements de groupes diff\'erentiels}
	Th\`ese de doctorat d'\'etat, Université de Provence, Marseille (1984)
	
	{{}
	\bibitem{Dud1968} Dudley, R.M., Distances of probability measures and random variables {\it Ann. Math. Stat.} {\bf 39 } 1563–1572 (1968)
}

	\bibitem{DN2007-1}
	Dugmore, D.; Ntumba, P.;On tangent cones of Fr\"olicher spaces
	\textit{Quaetiones mathematicae} \textbf{30} no1 (2007) 67-83
	
	\bibitem{DN2007-2} Dugmore, D.; Ntumba, P.; Cofibrations in the category of Fr\"olicher spaces: part I \textit{Homotopy, homology and applications} \textbf{9}, no 2 (2007) 413-444
	
	\bibitem{Ee}Eells, J.;  A setting for global analysis
	\textit{Bull. Amer. Math. Soc.} {\bf 72} 751-807 (1966)
	
	\bibitem{ERMR2017} Eslami Rad, A; Magnot, J-P.;  Reyes e. G.;  The Cauchy problem of the Kadomtsev-Petviashvili hierarchy with arbitrary coefficient algebra \textit{J. Nonlinear Math. Phys.} \textbf{24}  sup1.  103-120 (2017).
	
	\bibitem{FH2001} Fadell, E.R.; Husseini, S.Y.; \textit{Geometry and topology of configuration spaces} 
	Springer, Berlin (2001)
	

	\bibitem{FK} Fr\"olicher, A; Kriegl, A; \textit{Linear spaces and differentiation
		theory} Wiley series in Pure and Applied Mathematics, Wiley Interscience
	(1988)
	
	
	\bibitem{Ham}  Hamilton, R.S.; { The inverse function theorem of Nash and Moser}; 
	{\it Bull. Amer. Math Soc. (NS)} \textbf{7} (1984) 65-222 
	
	\bibitem{He1995} Hector, G.;
	G\'eom\'etrie et topologie des espaces diff\'eologiques
	, in: \textit{Analysis and Geometry in Foliated Manifolds
	(Santiago de Compostela, 1994)}, World Sci. Publishing,  55–80 (1995)
	
	\bibitem{HN1971} Hoghe-Nlend, H.; \textit{Th\'eorie des bornologies et applications} Lect. Notes in Math. \textbf{273} (1971)
	
	\bibitem{IgPhD} Iglesias-Zemmour, P.; \textit{Fibrations diff\'eologiques et homotopie},  PhD thesis, universit\'e de Provence (1985)
	
	\bibitem{Igdiff} Iglesias-Zemmour, P.; \textit{Diffeology} Mathematical Surveys and Monographs
	\textbf{185} (2013)
	
	\bibitem{IgK2014} Iglesias-Zemmour, P.; Karshon, Y.; Smooth Lie group actions are parametrized diffeological subgroups. \textit{Proc. AMS} {\bf 140} no 2,  731-739 (2012).
	
	\bibitem{KM} Kriegl, A.; Michor, P.W.; \textit{The convenient setting
		for global analysis} Math. surveys and monographs \textbf{53}, American
	Mathematical society, Providence, USA. (2000)
	
	
	
	\bibitem{Leandre2002} L\'eandre, R.; Analysis on Loop Spaces and Topology \textit{Mathematical Notes} \textbf{72} no 1, 212-229  (2002)
	
	\bibitem{dLS1}C. De Lellis, L. Székelyhidi Jr.
	The Euler equations as a differential inclusion.
	\textit{Ann. of Math. (2)} {\bf 170} no. 3, 1417–1436 (2009)
	
	\bibitem{dLS2}C. De Lellis, L. Székelyhidi Jr.
	On admissibility criteria for weak solutions of the Euler equations.
	{\it Arch. Ration. Mech. Anal.} {\bf 195}  no. 1, 225–260 (2010)
	\\
	and
	
	Errata to "On admissibility criteria for weak solutions of the Euler equations"
	
	\bibitem{dLS3} C. De Lellis, L. Székelyhidi Jr.
	The h-principle and the equations of fluid dynamics.
	{\it Bull. Amer. Math. Soc.} {\bf 49}, 347-375 (2012)
	
	\bibitem{Les} Leslie, J.; On a Diffeological Group Realization of
	certain Generalized symmetrizable Kac-Moody Lie Algebras \textit{J.
		Lie Theory} \textbf{13}, 427-442 (2003)
	

	
	
	\bibitem{Ma2006-3} Magnot, J-P.; Diff\'eologie sur le fibr\'e d'holonomie d'une connexion en dimension infinie
	\textit{C. R. 
		Math. Acad. Sci., Soc. R. Can.} 
	\textbf{28}, 
	no. 4, 121--127 (2006)
	
	
	
	\bibitem{Ma2013} Magnot, J-P.; Ambrose-Singer theorem on diffeological bundles and complete integrability of the KP equation; \textit{Int. J. Geom. Meth. Mod. Phys.} \textbf{ 10}, No. 9, Article ID 1350043, 31 p. (2013). 
	
	
	
	\bibitem{Ma2016-2} Magnot, J-P.; 
	Differentiation on spaces of triangulations and optimized triangulations.\textit{ 
	 5th international conference in mathematical modelling in physical science (IC-MSquare 2016) J. Phys: Conf. Ser.} \textbf{738} article ID 012088 (2016)
\bibitem{Ma2013-2}  Magnot, J-P.; The group of diffeomorphisms of a non-compact manifold is not regular
\textit{Demonstr. Math.} 51, No. 1, 8-16 (2018)
{{}
 \bibitem{Ma2020-1} Magnot, J-P.; On the domain of implicit functions in a projective limit setting without additional norm estimates. {\it Demonstr. Math.} {\bf 53}  no. 1, 112–120 (2020)}
\bibitem{MR2016} Magnot, J-P.; Reyes, E. G.; Well-posedness of the Kadomtsev-Petviashvili hierarchy, Mulase factorization and Frölicher Lie groups. 
{{} \textit{Ann. H. Poincar\'e} {\bf 21}, 1893--1945 (2020)
\bibitem{MW2017} Magnot, J-P.; Watts, J. A.; The diffeology of Milnor's classifying space. \textit{Top. Appl.} \textbf{232} 189-213 (2017)}
\bibitem{Nt2002} Ntumba, P.; DW Complexes and Their Underlying Topological Spaces
\textit{Quaestiones Math.}  {\bf 25},
119-134 (2002) 
	\bibitem{Olv} Olver, P.J.; \textit{Applications of Lie groups to differential equations (2nd edition)} GTM \textbf{107}, Springer (1993)
	
\bibitem{Om} Omori, H.; \textit{Infinite dimensional Lie groups} AMS translations of mathematical monographs \textbf{158} (1997)
	
	
	\bibitem{Pa} Palais, R.S.; Homotopy theory of infinite dimensional manifolds \textit{Topology} \textbf{5}  1-16 (1966)

	
	\bibitem{Pen}  Penot, J.-P.; {Sur le th\'eor\`eme de Frobenius }; {\it Bull. Soc. math. France } \textbf{98},
	47-80 (1970)
	
	\bibitem{pervova}
	Ekaterina Pervova, E.; Diffeological vector pseudo-bundles, \emph{Topology Appl.} \textbf{202} , 269--300 (2016).
	
	\bibitem{pervova2017}
	Pervova, Ekaterina
	Diffeological gluing of vector pseudo-bundles and pseudo-metrics on them. 
	\textit{Topology Appl.} \textbf{220}, 65-99 (2017). 
	\bibitem{pervova2}
	Pervova, E.;
	On the notion of scalar product for finite-dimensional diffeological vector spaces. 
	\textit{Electron. J. Linear Algebra} 
	\textbf{34}, 18-27 (2018). 
	
	\bibitem{pervova2018}  Pervova, E.; Diffeological Clifford algebras and pseudo-bundles of Clifford modules. To appear in \textit{Linear and Multilinear Algebra} 
	
	\bibitem{Prok} Prokhorov, Yu.V.; Convergence of random	processes and limit	theorems in	probability theory. \textit{Theor. Prob. Appl.} \textbf{1},
	157-214
	(1956)
	
	\bibitem{Rob} Robart, T.; Sur l'int\'egrabilit\'e des sous-alg\`ebres de
	Lie en dimension infinie; \textit{Can. J. Math.} \textbf{49} (4) (1997),
	820-839 
	
	
	\bibitem{Sch}  Scheffer, V.;  An inviscid flow with compact support in space-time. {\it J.	Geom. Anal.} {\bf 3} no 4 , 343–401 (1993)
	
	\bibitem{Sh1} Shnilerman, A.; On the nonuniqueness of weak solution of the Euler
	equation. \textit{Comm. Pure Appl. Math.} {\bf 50} no 12, 1261–1286 (1997)
	
	\bibitem{Sh2} Shnilerman, A.; Weak solutions with decreasing energy of incompressible
	Euler equations. \textit{Comm. Math. Phys.} \textbf{210} no 3, 541–603 (2000)
	
	\bibitem{Sou} Souriau, J.M.; Un algorithme g\'en\'erateur de structures quantiques; 
	\textit{Ast\'erisque}, Hors S\'erie, (1985) 341-399 
	
	
	\bibitem{Vil} Vilani, C.; Paradoxe de Scheffer-Shnirelman revu sous l’angle de l’intégration convexe (d’après C. De Lellis, L. Székelyhidi). {\it S\'eminaire Bourbaki} Exp. 1001 (November 2008)
	
	\bibitem{Vin2013}; Vinogradov, A.M. ; What are symmetries of nonlinear PDEs and what are they themselves? \texttt{ArXiv:1308.5861}
	
	\bibitem{Wa} Watts, J.; \textit{Diffeologies, differentiable spaces
		and symplectic geometry} PhD thesis, university of Toronto (2012) arXiv:1208.3634
\bibitem{We2017} Welker, K.; Suitable Spaces for Shape Optimization. \texttt{arXiv:1702.07579}	

\end{thebibliography}
\end{document}